\documentclass[leqno,12pt,twoside]{amsart}

\usepackage{times}

\usepackage{xy} % diagrams
\xyoption{all}
\usepackage{amsmath, amssymb, amsfonts, latexsym, graphics,mdwlist, enumitem}
\usepackage{subfig, floatflt}

\usepackage{a4}

\def\ch{\mathop{\mathrm{ch}}\nolimits}

\def\Coh{\mathop{\mathrm{Coh}}\nolimits}

\def\dim{\mathop{\mathrm{dim}}\nolimits}

\def\Ext{\mathop{\mathrm{Ext}}\nolimits}
\def\Hom{\mathop{\mathrm{Hom}}\nolimits}

\def\RlHom{\mathop{\mathbf{R}\mathcal Hom}\nolimits}
\def\RHom{\mathop{\mathbf{R}\mathrm{Hom}}\nolimits}

\def\num{\mathop{\mathrm{num}}\nolimits}

\def\Ob{\mathop{\mathrm{Ob}}}

\def\rk{\mathop{\mathrm{rk}}}

\def\supp{\mathop{\mathrm{supp}}}

\def\td{\mathop{\mathrm{td}}\nolimits}

\def\top{\mathop{\mathrm{top}}\nolimits}

\def\DT{\mathop{\mathrm{DT}}}
\def\PT{\mathop{\mathrm{PT}}}

\newenvironment{Prf}{\textit{Proof.}\/}{\hfill$\Box$}

\def\Stab{\mathop{\mathrm{Stab}}}
\def\StabP{\mathop{\mathrm{Stab}_{\mathrm{Pol}}}}
\def\into{\ensuremath{\hookrightarrow}}
\def\onto{\ensuremath{\twoheadrightarrow}}

\newcommand\TFILTB[3]{%
\xymatrix@=1pc{
{0 = {#1}_0} \ar[rr]&&
{{#1}_1} \ar[rr]\ar[ld] &&
{{#1}_2} \ar[r]\ar[ld] &
{\cdots} \ar[r] & { {#1}_{#3-1}} \ar[rr] &&
{{#1}_{#3} = {#1}} \ar[ld]
\\
& *{{#2}_1} \ar@{.>}[ul] &&
{{#2}_2} \ar@{.>}[ul] & &&&
{{#2}_{{#3}}} \ar@{.>}[ul]
}}

\def\abs#1{\lvert#1\rvert}

\newcommand\stv[2]{\left.\kern-\nulldelimiterspace
                \left\{#1\vphantom{#2}\,\right|#2\right\}}

\newtheorem{Def-s}[subsection]{Definition}
\newtheorem{Thm-s}[subsection]{Theorem}
\newtheorem{Prop-s}[subsection]{Proposition}
\newtheorem{Lem-s}[subsection]{Lemma}
\newtheorem{Thm}[subsubsection]{Theorem}
\newtheorem{Def}[subsubsection]{Definition}

\newtheorem{Prop}[subsubsection]{Proposition}

\newtheorem{Lem}[subsubsection]{Lemma}

\def\C{\ensuremath{\mathbb{C}}}
\def\D{\ensuremath{\mathbb{D}}}
\def\H{\ensuremath{\mathbb{H}}}

\def\Q{\ensuremath{\mathbb{Q}}}
\def\R{\ensuremath{\mathbb{R}}}
\def\Z{\ensuremath{\mathbb{Z}}}

\def\AA{\ensuremath{\mathcal A}}
\def\BB{\ensuremath{\mathcal B}}

\def\DD{\ensuremath{\mathcal D}}
\def\EE{\ensuremath{\mathcal E}}
\def\FF{\ensuremath{\mathcal F}}
\def\GG{\ensuremath{\mathcal G}}

\def\JJ{\ensuremath{\mathcal J}}
\def\KK{\ensuremath{\mathcal K}}
\def\LL{\ensuremath{\mathcal L}}
\def\MM{\ensuremath{\mathcal M}}
\def\NN{\ensuremath{\mathcal N}}
\def\OO{\ensuremath{\mathcal O}}
\def\PP{\ensuremath{\mathcal P}}

\def\QQ{\ensuremath{\mathcal Q}}
\def\TT{\ensuremath{\mathcal T}}

\def\ZZ{\ensuremath{\mathcal Z}}

\def\barp{\overline{p}}

\begin{document}
\title{Polynomial Bridgeland stability conditions and the large volume limit}

\author{Arend Bayer}
\address{University of Utah, Department of Mathematics, 155 South 1400 East,
Room 233, Salt Lake City, UT 84112}
\email{bayer@math.utah.edu}


\begin{abstract}
We introduce the notion of a polynomial stability condition, generalizing
Bridgeland stability conditions on triangulated categories.  We construct and
study a family of polynomial stability conditions for any normal projective
variety. This family includes both Simpson-stability, and large volume limits
of Bridgeland stability conditions.

We show that the PT/DT-correspondence relating stable pairs to
Donaldson-Thomas invariants (conjectured by Pandharipande and Thomas)
can be understood as a wall-crossing in our family of polynomial stability
conditions.  Similarly, we show that the relation between stable pairs and
invariants of one-dimensional torsion sheaves (proven recently by the same
authors) is a wall-crossing formula.
\end{abstract}

\maketitle

\setcounter{tocdepth}{1}
\tableofcontents

\section{Introduction}

In this article, we introduce polynomial stability conditions on triangulated
categories. They are a generalization of Bridgeland's notion of stability in
triangulated categories. The generalization is motivated by trying to
understand limits of Bridgeland's stability conditions; it allows for the
central charge to have values in complex polynomials rather than complex
numbers. 

While Bridgeland stability conditions have been constructed only in
dimension $\le 2$ and some special cases, we construct a family of
polynomial stability conditions on the
derived category of any normal projective variety. This family includes both
Simpson-stability of coherent sheaves, and stability conditions that 
we expect to be
the large volume limit of Bridgeland stability conditions.

We interpret both the PT/DT-correspondence conjectured 
in \cite{PT1}, and the relation between stable pair invariants and 
one-dimensional torsion sheaves proven in \cite{PT-BPS}, as a wall-crossing
phenomenon in our family of polynomial stability conditions.

\subsection{Bridgeland's stability conditions}
Since their introduction in \cite{Bridgeland:Stab}, stability conditions
for triangulated categories have drawn an increasing amount of interest
from various perspectives. They generalize the concept of stability 
from abelian categories to triangulated categories.

Originally, Bridgeland introduced the concept as an attempt to mathematically
understand Douglas' construction \cite{Douglas:stability} of $\Pi$-stability
of $D$-branes. Following Douglas' ideas, Bridgeland showed that the set of
stability conditions on $D^b(X)$ has a natural structure as a smooth manifold.
There are also various purely mathematical reasons to study the space of
stability conditions. 

\begin{Def}[\cite{Bridgeland:Stab}]	\label{def:Bridgeland-stability}
A stability condition on $D^b(X)$ is a pair $(Z, \PP)$ where
$Z \colon K(X) \cong K(D^b(X)) \to \C$ is a group homomorphism, and
$\PP$ is a collection of extension-closed subcategories
$\PP(\phi)$ for $\phi \in \R$, such that
\begin{enumerate}[label={(\alph*)}]
\item $\PP(\phi + 1) = \PP(\phi)[1]$,
\item $\Hom(\PP(\phi_1), \PP(\phi_2)) = 0$ for all 
$\phi_1 > \phi_2$,
\item if $0 \neq E \in \PP(\phi)$, then 
$Z(E) \in \R_{>0}\cdot e^{i \pi \phi}$, and
\item \label{cond:HN-filtration}
for every $0 \neq E \in D^b(X)$ there is a sequence
$\phi_1 > \phi_2 > \dots > \phi_n$ of real numbers and a sequence
of exact triangles
\[
\TFILTB E A n
\]
with $A_i \in \PP(\phi_i)$.
\end{enumerate}
\end{Def}
Objects of $\PP(\phi)$ are called semistable of phase $\phi$, and the
group homomorphism $Z$ is called the central charge.
We now restrict our attention to ``numerical'' stability conditions: these
are stability conditions for which $Z(E)$ is given
by numerical invariants of $E$, i.e. where $Z$ 
factors via the projection $K(D^b(X)) \to \NN(X) := \NN(D^b(X))$ 
to the numerical $K$-group.\footnote{The numerical $K$-group
$\NN(D^b(X))$ is the quotient
of $K(D^b(X))$ by the zero-space of the bilinear form
$\chi(E, F) = \chi(\RHom(E,F))$.}

\subsection{The space of stability conditions}
The role of $\PP$ (called ``slicing'') is easily understood, as it
naturally generalizes the notion of semistable objects in an abelian
category, together with the ordering of their slopes and the existence
of Harder-Narasimhan filtrations. The role of
$Z$ is less obvious; we will explain two aspects in the following
paragraphs.

It seems unsatisfactory that semistable objects in the derived category
have to be given explicitly, rather than characterized intrinsically by
a slope function. This deficiency is somewhat corrected by the following
observation:

Given a slicing $\PP$, consider the category
$\AA = \PP((0,1])$ generated by all semistable objects of phase
$0 < \phi \le 1$
and extensions. It can be seen that $\AA$ is the heart of a bounded
t-structure (and in particular an abelian category); the slicing is thus
a refinement of the datum of a bounded t-structure. Bridgeland showed that
this refinement is uniquely determined by $Z$:

\begin{Prop} \label{prop:Z-and-t}
\cite[Proposition 5.3]{Bridgeland:Stab}
To give a stability condition $(Z, \PP)$ is equivalent to giving the
heart $\AA \subset D^b(X)$ of a bounded t-structure, and a group
homomorphism $Z \colon K(\AA) \to \C$ with the following properties:
\begin{enumerate}[label={(\alph*)}]
\item \label{cond:positivity}
For every object $E \in \AA$, we have
$Z(E) \in \R_{>0} \cdot e^{i \pi \phi(E)}$ with
$0 < \phi(E) \le 1$.
\item
We say an object is $Z$-semistable if it has no subobjects
$A \into E$ with $\phi(A) > \phi(E)$. We require that
every object has a Harder-Narasimhan filtration with $Z$-semistable
filtration quotients.
\end{enumerate}
\end{Prop}
Given $\AA$ and $Z$, the semistable objects in the derived category
are the shifts of the $Z$-semistable objects in $\AA$. The positivity
condition \ref{cond:positivity} is somewhat delicate; for example, 
it can't be satisfied for the category of coherent sheaves on a
projective surface.

There is a natural topology on the space of slicings. However, only 
together with the central charge does the topological space of stability
conditions become a smooth manifold. One can paraphrase Bridgeland's
result as follows:
\emph{One can equip the space $\Stab(X)$ of
``locally finite''\footnote{\cite[Definition 5.7]{Bridgeland:Stab}}
numerical stability conditions with the
structure of a smooth manifold, such that the forgetful map $\ZZ \colon
\Stab(X) \to \NN(X)^*$,  $(Z, \PP) \mapsto Z$ gives local coordinates at every
point.}
In other words, a stability condition can be deformed by deforming its
central charge. 

The space $\Stab(X)$ is closely related to the moduli space of $N=2$
superconformal field theories, see \cite{Bridgeland:spaces}.
The existence of $\ZZ$ has interesting implications on the group of
auto-equivalences of $D^b(X)$, as one can study its induced action
on $\ZZ$, see e.g. \cite{Bridgeland:K3} and \cite{HMS:generic_K3s}.

\subsection{Reconstruction of $X$ from $D^b(X)$}
\label{sect:reconstruction}
If the canonical bundle $\omega_X$, or its inverse $\omega_X^{-1}$, of a smooth
variety $X$ is ample, then the variety can be reconstructed from its bounded
derived category, see \cite{Bondal-Orlov}.
Without the assumption of ampleness, this statement is
wrong, and the proof already breaks down at its first step: the intrinsic
characterization of point-like objects in $D^b(X)$ (the shifts $\OO_x[j]$ of
skyscraper sheaves for closed points $x \in X$) by the action of the
Serre-functor.

However, the mathematical translation of ideas by Aspinwall, originally
suggested in \cite{Aspinwall:points}, suggests that a stability condition
provides exactly the missing data to characterize the point-like objects.
Inside the space of stability conditions, there should be a special chamber,
which we will call the ample chamber, with the following property: When $(Z,
\PP)$ is a stability condition in the ample chamber, and $E \in D^b(X)$ an
object with class $[E] = [\OO_x]$ in the numerical $K$-group, then $E$ is
$(Z, \PP)$-stable if and only if $E$ is isomorphic to the shift of a
skyscraper sheaf $[\OO_x]$.  One could then reconstruct $X$ as the moduli
space of $(Z, \PP)$-stable objects.

Moving to a chamber of the space of stability conditions adjacent to the ample
chamber, the moduli space $\widetilde X$ of semistable objects of the same
class $[\OO_x]$ comes with a fully faithful functor
$\Phi \colon D^b(\widetilde X) \to D^b(X)$ induced
by the universal family. This suggests that $\widetilde X$ could be a
birational model of $X$ with isomorphic derived category (e.g. a flop),
it could be isomorphic to $X$ with $\Phi$ being a non-trivial
auto-equivalence of $X$, or it could be a birational contraction or
a flip of $X$. It seems an intriguing question to what extent the birational
geometry of $X$ can be captured by this phenomenon.

This suggestion is consistent with many of the known examples of Bridgeland
stability conditions.  Maybe most convincing is the case of a crepant resolution
$Y \to \C^3/G$ of a three-dimensional Gorenstein quotient singularity.  The
results in \cite{Alastair-Ishii} can be reinterpreted as saying that every
other crepant resolution $Y' \to \C^3/G$ can be constructed as a moduli
space of Bridgeland-stable objects in $D^b(Y)$; see also
\cite{Toda:stab-crepant_res} for the local construction of a flop
along these lines.

\subsection{Examples of stability conditions}
The existence of stability conditions on $D^b(X)$ for $X$ a smooth,
projective variety has only been shown in very few cases:
\begin{itemize*}
\item
For a smooth curve $C$, stability conditions have been constructed
in \cite{Bridgeland:Stab}, and $\Stab(C)$ has been described by
\cite{Macri:curves, Okada:P1}; in
\cite{BK:elliptic} the case of singular curves of genus one was
considered.
\item For the case of a K3 surface $S$, Bridgeland completely described one
connected component of $\Stab(S)$ in \cite{Bridgeland:K3} (including a complete
description of the ample chamber). In \cite{MMS:inducing},
the authors study the space of stability conditions on Kummer and Enriques
surfaces.  For arbitrary smooth projective surfaces, stability conditions
have recently been constructed in \cite{Aaron-Daniele}.
\item If $D^b(X)$ has a complete exceptional collection, then 
stability conditions exist by \cite{Macri:stability-examples}.
\end{itemize*}
For complex non-projective tori, stability conditions have been studied
in \cite{Sven:generic_tori}.

\subsection{Stability conditions related to $\sigma$-models}
Let $X$ be a smooth projective variety.
Following ideas in the physical literature (see
\cite{Douglas:stability, Aspinwall-Douglas:stability, Aspinwall:points,
Aspinwall-Lawrence:DC-zero-brane}), it should be possible to construct
stability conditions on $D^b(X)$ coming from the non-linear $\sigma$-model 
associated to $X$. At least for an open subset of these stability conditions,
skyscraper sheaves of points should be stable. Further, it is known
how the central charge should depend on the complexified K\"ahler
moduli space: if $\beta \in H^2(X)$ is an arbitrary class, and
$\omega \in H^2(X)$ and ample class, then the central charge
should be given as
\begin{equation} 		\label{eq:central-charge}
Z_{\beta, \omega}(E) = -\int_X e^{-\beta - i \omega} \cdot \ch(E) \sqrt{\td X}.
\end{equation}
However, in general not even a matching t-structure whose heart $\AA$
would satisfy the positivity condition \ref{cond:positivity} of
Proposition \ref{prop:Z-and-t} is known; in fact no example of
a stability condition on a projective Calabi-Yau threefold is known.

\subsection{Polynomial stability conditions} \label{sect:intro-pol-stab}
However, if we replace $\omega$ by $m \omega$ and let
$m \to +\infty$ (this is the large volume limit), then a matching t-structure
can be constructed: If $E$ is a coherent sheaf and $d$ its dimension of
support, then $Z(E)(m) \to -(-i)^d \cdot \infty$ as $m \to \infty$. Thus the
central charge $Z(E[\lfloor \frac d2\rfloor])(m)$ of the shift of $E$ will go
to $-\infty$ or $i \infty$; this suggests that a t-structure can be
constructed by a filtration of dimension of support, i.~e.~a t-structure of
perverse coherent sheaves. However, the limit of the phase $\phi(E)(m)$
is too coarse as an information to characterize semistable objects;
instead, it is more natural to consider the central charge $Z_{\beta,
m\omega}$ given by equation (\ref{eq:central-charge}) as a polynomial in $m$:
then we can say a perverse coherent sheaf $E$ is semistable if there is no
perverse coherent subsheaf $A \into E$ with $\phi(A)(m) > \phi(E)(m)$
for $m$ being large.

Motivated by this observation, we introduce a notion of polynomial stability
condition in definition \ref{def:polstab}.  It allows the central charge to
have values in polynomials $\C[m]$ instead of $\C$; accordingly, the slicing
$\PP$ has to depend not on real numbers, but on phases of polynomials
(considered for $m \gg 0$).  It gives a precise meaning to the notion of a
``stability condition in the limit of $m \to \infty$''.

\subsection{Results}
Our main result is Theorem \ref{mainthm}. It shows the existence of a 
family of polynomial stability conditions for every normal projective
variety. Its associated bounded t-structure is a t-structure of
perverse coherent sheaves.
The family contains stability conditions corresponding
to Simpson stability (see section \ref{sect:Simpson}), and stability
conditions that should be the large volume limit of Bridgeland
stability conditions (see section \ref{sect:large-volume}).

In the case of surfaces, Proposition \ref{prop:large-volume} makes the last
statement precise: the polynomial stability condition $(Z, \PP)$ at the large
volume limit is the limit of Bridgeland stability conditions $(Z_m, \PP_m)$,
depending on $m$, in the sense that objects are $\PP$-stable if and only if
they are $\PP_m$-stable for $m \gg 0$, and the Harder-Narasimhan filtration
with respect to $\PP$ is the same as the Harder-Narasimhan filtration with
respect to $\PP_m$ for $m \gg 0$.

The polynomial stability conditions provide many new t-structures on the
derived category of a projective variety.\footnote{The t-structures used
in the construction are those described in \cite{Bezrukavnikov:perverse},
but tilting with respect to different phase functions yields new torsion
pairs, and thus new t-structures, in the same way that Gieseker- or
slope-stability yield new t-structures by tilting the category
of coherent sheaves.}
They might help to construct
Bridgeland stability conditions on higher-dimensional varieties.

With Proposition \ref{prop:X-as-modspace}, we observe that the
polynomial stability conditions constructed in Theorem \ref{mainthm}
are ``ample'' in the sense of section \ref{sect:reconstruction}:
$X$ can be reconstructed from $D^b(X)$, the stability condition, and
the class of $[\OO_x] \in \NN(X)$ as a moduli space of semistable objects.

\subsection{PT/DT-correspondence as a wall-crossing}
In \cite{PT1}, the authors introduced new invariants of stable pairs
on smooth projective threefolds. In the Calabi-Yau case, they conjecture
a simple relation between their generating function and the generating
function of Donaldson-Thomas invariants (introduced in \cite{MNOP1}).
With Proposition \ref{prop:PTDT}, we show that this relation can be understood
as a wall-crossing phenomenon (in the sense of \cite{Joyce4}) in a family of
polynomial stability conditions.

Similarly, we show in section \ref{sect:PT-BPS} that the relation between
stable pair invariants and invariants counting one-dimensional torsion
sheaves can be understood as a wall-crossing formula.

\subsection{The space of polynomial stability conditions}
In section \ref{sect:space}, we discuss to what extent the deformation
result by Bridgeland carries over to our situation. We introduce
a natural topology on the set of polynomial stability conditions,
and show that the forgetful map
\[
\ZZ \colon \Stab\nolimits_{\mathrm{Pol}}(X) \to \Hom(\NN(X), \C[m]),
\quad (Z, \PP) \mapsto Z
\]
is continuous and locally injective. Under a strong 
local finiteness assumption, we can also show that it
is a a local homeomorphism.

\subsection{Notation}

If $\Sigma$ is a set of objects in a triangulated category $\DD$
(resp. a set of subcategories of $\DD$), we write
$\langle \Sigma \rangle$ for the full subcategory generated by $\Sigma$ and
extensions; i.e. the smallest full subcategory of $\DD$ that is closed
under extensions and contains $\Sigma$ (resp. contains all subcategories
in $\Sigma$).

We will write $\H \subset \C$ for the semi-closed upper half plane
\[ \H = \stv{z \in \C}{z \in \R_{>0} \cdot e^{i \pi \phi(z)}, 0 < \phi(z)
\le 1},\]
and $\phi(z)$ for the phase of $z \in \H$.

\subsection{Acknowledgments}
I would like to thank
Yuri I. Manin for originally suggesting the viewpoint of section
\ref{sect:reconstruction}, Richard Thomas for discussions related to section
\ref{sect:PTDT}, and Aaron Bertram, Nikolai Dourov, Daniel Huybrechts, Yunfeng
Jiang, Davesh Maulik and Gueorgui Todorov for useful comments and discussions.

Some of the stability conditions constructed in this article have also
been constructed independently by Yukinobu Toda in \cite{Toda:limit-stable},
namely the stability conditions at the large volume limit of Calabi-Yau
threefolds. In particular, Toda also explains the key formula
(\ref{eq:BPS-relation}) as wall-crossing formula in his family of stability
conditions. The complete family of stability conditions considered by Toda
lives on a wall of the space of polynomial stability conditions considered
here.

\section{Polynomial stability conditions}

\subsection{Example: Simpson/Rudakov stability as a polynomial stability
condition}
\label{sect:Simpson}

Before giving the precise definition of polynomial stability conditions,
we give an example that is more easily constructed than
the large volume limit considered in the introduction, which will
hopefully motivate the definition.

\begin{figure}[htb]
\begin{center}
        \includegraphics{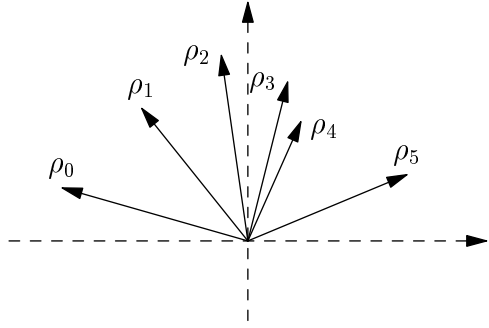}
\caption{Stability vector for Simpson stability}
\label{fig:simpson}
\end{center} \end{figure}

Let $\AA = \Coh X \subset D(X)$ be the standard heart in the derived
category of a projective variety
$X$ with a chosen ample line bundle $\LL$. Pick complex numbers
$\rho_0, \rho_1, \dots,
\rho_n$ in the open upper half plane $H$
with $\phi(\rho_0) >
\phi(\rho_1) > \dots > \phi(\rho_n)$ as in figure \ref{fig:simpson}.
For any coherent sheaf $E \in \Coh X$,
let $\chi_E(m) = \sum_{i=0}^n a_i(E) m^i$ be the Hilbert polynomial with
respect to $\LL$. We define the central charge by
\[
Z(E)(m) = \sum_{i=0}^n \rho_i a_i(E) m^i.
\]
Then $Z(E)(m) \in \H$ for $E$ nontrivial and $m \gg 0$, and we can 
consider the phase $\phi(E)(m) \in (0, 1]$. We say that a sheaf
if $Z$-stable if for every subsheaf $A \into E$, we have
$\phi(A)(m) \le \phi(E)(m)$ for $m \gg 0$.

Then a sheaf $E$ is $Z$-stable if and only if it is a Simpson-stable
sheaf; this is most easily seen by using Rudakov's reformulation in
\cite{Rudakov:stability}. In particular, stability does not depend on the
particular choice of the $\rho_i$.
In order not to lose any information, we should consider the phase
of a stable object 
to be the function $\phi(E)(m)$ defined for $m \gg 0$
rather than the limit $\lim_{m \to \infty} \phi(E)(m)$; in other
words, we consider its phase to be the function germ 
\[ \phi(E) \colon (\R \cup \{+\infty\}) \to \R. \]
Then we can define an object $E \in D^b(X)$ to be stable if and only 
if it is isomorphic to the shift $F[n]$ of $Z$-stable sheaf;
its phase is given by the function germ $\phi(F) + n$.

Combining the Harder-Narasimhan filtrations of arbitrary sheaves
with respect to Simpson stability with the filtration of a complex
by its cohomology sheaves, we obtain a filtration of an arbitrary
complex similar to the filtration in part \ref{cond:HN-filtration} of 
definition \ref{def:Bridgeland-stability}.

\subsection{Slicings}

\begin{Def}						\label{def:slicing}
Let $(S, \succeq)$ be a linearly ordered set, equipped with an
order-preserving bijection $S \to S, \phi \mapsto \phi + 1$ (called the shift)
satisfying $\phi + 1 \succeq \phi$.
An $S$-valued slicing of a triangulated category $\DD$ is given by full
additive extension-closed subcategories $\PP(\phi)$
for all $\phi \in S$, such that the following properties are satisfied:
\begin{enumerate}[label={(\alph*)}]
\item
For all $\phi \in S$, we have 
$\PP(\phi + 1) = \PP(\phi)[1]$.
\item \label{cond:Homvanishing}
If $\phi \succ \psi$ for $\phi, \psi \in S$, and $A \in \PP(\phi), B \in
\PP(\psi)$, then
$\Hom (A, B) = 0$
\item \label{cond:slicing-HN-filtration}
For all non-zero objects $E \in \DD$, there is a finite sequence
$\phi_1 \succ \phi_2 \succ \dots \succ \phi_n$ of elements in $S$, and a
sequence of exact triangles
\begin{equation}					\label{HN-system}
\TFILTB E A n
\end{equation}
\end{enumerate}
with $A_i \in \PP(\phi_i)$.
\end{Def}
This was called 
``stability data'' or ``t-stability'' in \cite{Christmas}.
If $S = \Z$, this notion is equivalent to a bounded t-structure (see
\cite[Lemma 3.1]{Bridgeland:K3}), and for $S = \R$, it is a 
``slicing'' as defined in \cite{Bridgeland:Stab}.
The objects
in $\PP(\phi)$ are called semistable of phase $\phi$. 
The sequence of exact triangles in part \ref{cond:slicing-HN-filtration}
is also called the
Harder-Narasimhan filtration of $E$. If a Harder-Narasimhan filtration
exists, then condition \ref{cond:Homvanishing} forces it to be unique.

\begin{Def}
The set $S$ of polynomial phase functions is the set of continuous
function germs
\[ \phi \colon (\R \cup \{+\infty\}, +\infty) \to \R \]
 such that there exists a polynomial $Z(m) \in \C[m]$ 
with $Z(m) \in \R_{> 0} \cdot e^{\pi i \phi(m)}$ for $m \gg 0$. 
It is linearly ordered by setting
\begin{eqnarray*}
\phi \prec \psi &\Leftrightarrow& \phi(m) < \psi(m) \quad
\text{for $0 \ll m < +\infty$},
\end{eqnarray*}
and its shift $\phi \mapsto \phi + 1$ is given by point-wise addition.
\end{Def}
The condition that $\phi, \psi$ can be written as arguments of
polynomial functions guarantees that either $\phi \succ \psi$ 
or $\phi \preceq \psi$ holds; given
$Z(m)$, the function $\phi(m)$ is of course determined up to an
even integer constant.

From now on, $S$ will be the set of polynomial phase functions.
In our construction, $S$-valued slicings will play the role of
$\R$-valued slicings in Bridgeland's construction.

The following easy lemma is implicitly used in both \cite{Bridgeland:Stab} and
\cite{Christmas}, but we will make it explicit:
\begin{Lem}
Let $S_1, S_2$ be two linearly ordered sets equipped with shifts
$\tau_1, \tau_2$,
and let $\pi \colon S_1 \to S_2$ be a morphism of ordered sets commuting
with $\tau_1, \tau_2$. Then $\pi$
induces a push-forward
of stability conditions as follows: If $\PP$ is an $S_1$-valued slicing, then
$\pi_*\PP(\phi_2)$ for some $\phi_2 \in S_2$ is defined as
$\langle \left\{\PP(\phi_1) \mid \pi(\phi_1) = \phi_2\right\} \rangle$.
\end{Lem}
The proof is an exercise in the use of the octahedral axiom.

We will make use of the following push-forwards: By the
projection $\pi \colon S \to \R, \phi \mapsto \phi(\infty)$, we obtain an
$\R$-valued slicing from
every $S$-valued slicing. Further, for each $\phi_0 \in S$ we get a
projection $\pi^{\phi_0} \colon S \to \Z,
\phi \mapsto \max_{n \in \Z} \phi_0 + n \preceq \phi$
(we could also choose $\phi \mapsto \max_{n \in \Z} \phi_0 + n \prec \phi$).
This produces a bounded t-structure from every $S$-valued slicing; in
other words, an $S$-valued slicing is a refinement of a bounded
t-structure, breaking up the category into even smaller slices.

For any interval $I$ in the set of phases, we get an
extension-closed subcategory
$\PP(I) = \langle \left\{\PP(\phi) \mid \phi \in I\right\} \rangle$.
In the case of an $S$-valued slicing, the categories
$\PP([\phi, \phi+1))$ and $\PP((\phi, \phi+1])$ are abelian, as they
are the hearts of the t-structures constructed in the last paragraph.
The proof for these statements carries over literally from the one
given by Bridgeland: we can include these categories into the
abelian category $\PP([\phi, \phi+1))$. The slices $\PP(\phi)$ are abelian.

\subsection{Central charge}
We now come to the main definition:
\begin{Def}				\label{def:polstab}
A polynomial stability condition on a triangulated category $\DD$ is given
by a pair $(Z, \PP)$, where $\PP$ is an $S$-valued slicing
of $\DD$, and $Z$ is a group homomorphism
$Z \colon K(\DD) \to \C[m]$, with the following property: 
if $0 \neq E \in \PP(\phi)$, then
\[ Z(E)(m) \in \R_{>0} \cdot e^{\pi i \phi(m)} \]
for $m \gg 0$.
\end{Def}
In the case where $Z$ maps to constant polynomials $\C \subset \C[m]$,
this is equivalent to Bridgeland's notion of a stability condition. Similarly
to that case, a polynomial stability condition can be constructed from
a bounded t-structure and a compatible central charge $Z$:

\begin{Def}					\label{def:centeredslope}
A polynomial stability function on an abelian category $\AA$
is a group homomorphism $Z \colon K(\AA) \to \C[m]$ such that there exists
a polynomial phase function $\phi_0 \in S$ with the following property:

For any $0 \neq E \in \AA$, there is a polynomial phase function
$\phi(E)$ with 
$\phi_0 \prec \phi(E) \preceq \phi_0 +1$
and $Z(E)(m) \in \R_{>0} \cdot e^{\pi i \phi_E(m)}$ for $m \gg 0$.
\end{Def}
This definition allows slightly bigger freedom than requiring
$Z(E)(m) \in \H$ for $m \gg 0$.

We call $\phi(E) \in S$ the phase of $E$; the function
$\Ob \AA \setminus \{0\} \to S$,
$E \mapsto \phi(E)$ is a slope function in the sense that it satisfies
the see-saw property on short exact sequences (cf. \cite{Rudakov:stability}).
An object $0 \neq E$ is called
\emph{semistable with respect to $Z$} if
for all subobjects $0 \neq A \subset E$, we have
$\phi(A) \preceq \phi(E)$; equivalently, if for every
quotient $E \onto B$ in $\AA$ we have $\phi(E) \preceq \phi(B)$.
We say that a stability function has the Harder-Narasimhan property
if for all $E \in \AA$, there is a finite filtration
$0 = E_0 \into E_1 \into \dots \into E_n = E$ such that
$E_i/E_{i-1}$ are semistable with slopes
$\phi(E_1/E_0) \succ \phi(E_2/E_1) \succ \dots \succ \phi(E_n/E_{n-1})$.

Finally, note that the set of polynomials $Z(E)$ for which a polynomial
phase function $\phi(E)$ as in the above definition exist forms a convex cone
in $\C[m]$. Its only extremal ray is the set of polynomials with
$\phi(E) = \phi_0 + 1$. This is an important reason why many of the proofs of
\cite{Bridgeland:Stab} carry over to our situation.

We restate two propositions by Bridgeland in our context; the proofs
are identical to the ones given by Bridgeland:
\begin{Prop} \cite[Proposition 5.3]{Bridgeland:Stab}
\label{prop:centeredslope} 
Giving a polynomial stability condition on $\DD$ is equivalent to giving a
bounded t-structure on $\DD$ and a polynomial stability function
on its heart with the Harder-Narasimhan property.
\end{Prop}

The following proposition shows that the Harder-Narasimhan property
can be deduced from finiteness assumption of $\AA$ with respect to
$Z$:
\begin{Prop}  \cite[Proposition 2.4]{Bridgeland:Stab}
					\label{prop:HNproperty}
Assume that $\AA$ is an abelian category, $Z \colon K(\AA) \to \C[m]$
a polynomial stability function, and that they satisfy the following
chain conditions:
\begin{description*}
\item[$Z$-Artinian]
There are no infinite chains of subobjects
\[ \dots \into E_{j+1} \into E_j \into \dots \into E_2 \into E_1 \]
with $\phi(E_{j+1}) \succ \phi(E_j)$ for all $j$.
\item[$Z$-Noetherian] 
There are no infinite chains of quotients
\[ E_1 \onto E_2 \onto \dots \onto E_j \onto E_{j+1} \onto \dots  \]
with $\phi(E_{j}) \succ \phi(E_{j+1})$ for all $j$.
\end{description*}
Then $\AA, Z$ have the Harder-Narasimhan property.
\end{Prop}

\section{The standard family of polynomial stability condition}
\label{sect:standard-stability}

In this section, we will construct a standard family of stability conditions
on the bounded derived category $D^b(X)$ of an arbitrary normal projective
variety $X$. Let $n$ be the dimension of $X$.

\subsection{Perverse coherent sheaves}
The t-structures relevant for our stability conditions are t-structures
of perverse coherent sheaves.  The theory of perverse coherent sheaves is
apparently originally due to Deligne, and has been developed by Bezrukavnikov
\cite{Bezrukavnikov:perverse} and Kashiwara \cite{Kashiwara}.
We will need only a special case of perverse coherent sheaves, which are
given by filtrations of dimension.

\begin{Def}					\label{def:perversity}
A function $p \colon \{0, 1, \dots, n\} \to \Z$ is called
a perversity function if $p$ is monotone decreasing, and if
$\barp \colon \{0, 1, \dots, n\} \to \Z$ (called the dual perversity)
given by $\barp(d) = -d - p(d)$ is also monotone decreasing.
\end{Def}
In other words we require that $p(d) \ge p(d+1) \ge p(d) - 1$.  Given a
perversity function in the above sense, the function
$X^{\top} \to \Z, x \mapsto p(\dim x)$ is a \emph{monotone and comonotone}
perversity function in the sense of \cite{Bezrukavnikov:perverse}.

Let $\AA^{p, \le k}$ be the following increasing
filtration of $\Coh X$ by abelian subcategories:
\[
\AA^{p, \le k} = \stv{\FF \in \Coh X}{p(\dim \supp \FF) \ge -k}
\]
\begin{Thm}\cite{Bezrukavnikov:perverse, Kashiwara}					 \label{thm:tstruct}
If $p$ is a perversity function, then the following pair defines a 
bounded t-structure on $D^b(X)$:
\begin{align}
D^{p, \le 0}
  & = \stv{E \in D^b(X)}{H^{-k}(E) \in \AA^{p, \le k}
				\quad \text{for all $k \in \Z$}}
					\label{eq:le0} \\
D^{p, \ge 0} & = \left\{ E \in D^b(X) \mid \Hom(A, E) = 0
	\quad \text{for all $k \in \Z$ and $A \in \AA^{p, \le k}[k+1]$} \right\}
					\label{eq:ge0}
\end{align}
\end{Thm}
This description is slightly different to the one given in
\cite{Bezrukavnikov:perverse, Kashiwara} but easily seen to be equivalent.
Once $D^{p, \le 0}$ is given, $D^{p, \ge 0}$ is of course
determined as the right-orthogonal complement of $D^{p, \le -1}$.
Our notation is somewhat intuitive as $\AA^{p, \le k}$ can be
recovered as $\AA \cap D^{p, \le k}$, which completely determines
the t-structure.

Objects in the heart $\AA^p = D^{p, \ge 0} \cap D^{p, \le 0}$ are called
perverse coherent sheaves.

\subsection{Construction of polynomial stability conditions}
\label{sect:mainthm}

\begin{figure}[htb]
\begin{center}
        \includegraphics{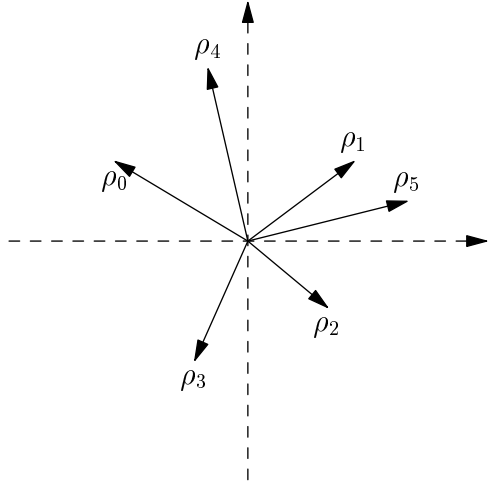}
\caption{A stability vector with associated perversity function
$p(0) = p(1) = 0$, $p(2) = p(3) = -1$, $p(4) = p(5)=-2$}
\label{fig:polstab}
\end{center} \end{figure}

\begin{Def}
A stability vector $\rho$ is a sequence
$(\rho_0, \rho_1, \dots, \rho_n) \in (\C^*)^{n+1}$
of non-zero complex numbers such that
$\frac {\rho_d}{\rho_{d+1}}$ is in the open upper half plane for $0 \le d \le n-1$.

Given a stability vector $\rho$, we call
$p \colon \{0, 1, \dots, n\} \to \Z$ a perversity
function associated to $\rho$ if it is a perversity function
satisfying
$(-1)^{p(d)} \rho_d  \in \H$ for all $0 \le d \le n$.
\end{Def}
Such $p$ is uniquely determined by $p(0)$, and given $p(0)$ such
a perversity function exists if $p(0)$ is of the correct parity;
see figure \ref{fig:polstab} for an example on a 5-fold. The
number $p(0)- p(d)$ counts how often the piecewise linear path
$\rho_0 \to \rho_1 \to \dots \to \rho_d$ crosses the real line.
We will construct stability conditions by giving a polynomial
stability functions on $\AA^{p}$. 

In the following, a Weil divisor $\omega \in A^1(X)_\R$ is called ample
if for any effective class $\alpha \in A_d(X)$, we have
$\omega^d \cdot \alpha > 0$.

\begin{Thm}					\label{mainthm}
Let the data $\Omega = (\omega, \rho, p, U)$ be given, consisting of
\begin{itemize*}
\item
an ample class $\omega \in A^1(X)_\R$,
\item a stability vector $\rho = (\rho_0, \dots, \rho_n)$,
\item a perversity function $p$ associated to $\rho$, and
\item  a unipotent operator $U \in A^*(X)_\C$ 
(i.e. $U = 1 + N$ where $N$ is concentrated in positive degrees).
\end{itemize*}
Let $Z_\Omega \colon K(X) \to \C[m]$ be the following central charge:
\[ Z_\Omega(E)(m) = \int_X \sum_{d=0}^n \rho_d \omega^d m^d \cdot \ch(E) \cdot U
\]
Then $Z_\Omega(E)(m)$ is a polynomial stability function for $\AA^p$
with the Harder-Narasimhan property. 
\end{Thm}
By Proposition \ref{prop:centeredslope}, this gives a polynomial stability
condition $(Z_\Omega, \PP_\Omega)$ on $D^b(X)$.

We will drop the subscript $\Omega$ from the notation.
In this section we will just 
prove that $Z$ is a polynomial stability function according
to definition \ref{def:centeredslope} with $\phi_0 = \epsilon$ for
some small constant $\epsilon \ge 0$. In other words,
we have to prove that for every $E \in \AA^p$, we have
$Z(E)(m) \in e^{i \epsilon} \cdot \H$ for $m \gg 0$.
The proof of the existence of Harder-Narasimhan filtrations will be
postponed until section \ref{sect:HN-proof}.

We start the proof with the following immediate observation:
\begin{Lem} 					\label{lem:dimsupport}
Given a non-zero object $E \in \AA^p$, let
$k$ be the largest integer such that $H^{-k}(E) \neq 0$, and let
$d$ be the dimension of support of $H^{-k}(E)$.
Then $p(d) = -k$, the sheaf $H^{-k}(E)$ has no torsion in dimension $d'$
whenever $p(d') > -k$, and all other cohomology sheaves of $E$ are supported
in smaller dimension.
\end{Lem}
We call $d$ the dimension of support of $E$.

\begin{Prf}
By $E \in D^{p, \le 0}$ we have 
$ p \left(\dim \supp H^{-k(E)}(E) \right) \ge -k$.
The claim follows from $E \in D^{p, \ge 0}$ and
\[ \Hom(\AA^{p, \le k-1}, H^{-k}(E))
= \Hom(\AA^{p, \le k-1}[k],  H^{-k}(E)[k])
= \Hom(\AA^{p, \le k-1}[k],  E) = 0 \]
\end{Prf}

Choose $\epsilon > 0$ such that $(-1)^{p(d)} \rho_d$ is in the interior
of $\H_\epsilon = e^{i \epsilon} \cdot \H$ for all $d$; we will first show that
$Z(E)(m) \in \H_\epsilon$ for $m \gg 0$.

Let $k$ be as in the lemma, and $d = \dim \supp H^{-k}(E)$. Since
all other cohomology sheaves of $E$ are supported in lower dimension,
we have
\[ \left(\ch(E) \cdot U \right)_{n-d} = (-1)^k \ch_{n-d}(H^{-k}(E)).\]
Since $\omega$ is ample and $\ch_{n-d}(H^{-k}(E))$ is effective,
the intersection product
$a := \int_X \omega^d \cdot \ch_{n-d}(H^{-k}(E))$ is positive. Thus the
leading term of $Z(E)(m)$ is $ a (-1)^d \rho_d m^d $.  Since
$a (-1)^d \rho_d \in \H_\epsilon$, the same must hold for $Z(E)(m)$
and large $m$.

\subsection{Dual stability condition}
					\label{sect:dual-stability1}

Let $\omega_X$ be a local dualizing complex of $X$, and let
\[ \D \colon D^b(X) \to D(X), \quad E \mapsto \RlHom(E, \omega_X) \]
be the associated dualizing functor. Let $D$ be such that
$\omega_X|_{X^{\text{smooth}}}$ is the shift of a line bundle 
by $D$.

To every polynomial stability condition $(Z_\Omega, \PP_\Omega)$ of Theorem 
\ref{mainthm} one can explicitly construct a stability
condition dual to $(Z, \PP)$ under $\D$. In the case where $X$ is not smooth,
this will be a stability condition on $\D(D^b(X))$ rather than $D^b(X)$; 
however, its associated heart is still given by a category of perverse
coherent sheaves as described earlier.

Let $P \colon A_*(X) \to A_*(X)$ be the parity operator acting by
$(-1)^{n-d}$ on $A_d(X)$.
Given the data $\Omega = (\omega, \rho, p, U)$ as in
Theorem \ref{mainthm}, we define
the dual data $\Omega^* = (\omega, \rho^*, \barp, U^*)$ by
$\rho^*_d = (-1)^{D+d} \overline{\rho_d}$,
$U^* = (-1)^D\ch(\omega_X)^{-1} \cdot P(\overline{U})$.
Consider the central charge $Z_{\Omega^*} \colon K(X) \to \C[m]$ defined
by the same formula as $Z_\Omega$ in \ref{mainthm}.

\begin{Prop} 					\label{prop:dual-stability}
The central charge $Z_{\Omega^*}$ induces a polynomial stability function
on $\D(\AA^p)$. 
The induced polynomial stability condition $(Z_{\Omega^*}, \PP_{\Omega^*})$
is dual to $(Z_\Omega, \PP_\Omega)$ in the following sense:
\begin{enumerate}[label={(\alph*)}]
\item An object $E$ is $(Z_\Omega, \PP_\Omega)$-stable if and only if $\D(E)$ is
$(Z_{\Omega^*}, \PP_{\Omega^*})$-stable.
\item If $E, F$ are $(Z_\Omega, \PP_\Omega)$-stable, then
\[ \phi(E) \succ \phi(F) \Leftrightarrow \phi(\D(E)) \prec \phi(\D(F)) \]
\item The Harder-Narasimhan filtration of $\D(E)$ with respect to
$(Z_{\Omega^*}, \PP_{\Omega^*})$ is obtained from that of $E$ with respect to 
$(Z_\Omega, \PP_\Omega)$ by dualization.
\end{enumerate}
\end{Prop}
By the uniqueness of HN filtrations, (a) and (b) imply (c). 
The proof of (a) and (b) will also be postponed until section
\ref{sect:HN-proof}.

\section{The large volume limit}
\label{sect:large-volume}

Fix $\beta \in A^1(X)_\R$ and an ample class $\omega_0 \in A^1(X)_\R$.
Let $\rho_d = - \frac{(-i)^d}{d!}$ and let $U = e^{-\beta}\cdot \sqrt{\td(X)}$.
Then
$p(d) = -\lfloor \frac d2 \rfloor$ is a perversity function associated
to $\rho = (\rho_0, \dots, \rho_n)$, and
the central charge $Z = Z_{\Omega}$ of Theorem \ref{mainthm} for
$\Omega = (\omega, \rho, p, U)$ is given by
\[
Z(E)(m) =  - \int_X e^{-\beta - i m\omega} \cdot \ch(E) \sqrt{\td(X)} \]\[
\]
This is the central charge $Z_{\beta, m \omega}$ discussed in section
\ref{sect:intro-pol-stab} as the central charge at the large-volume limit. 

This stability condition has many of the properties predicted by physicists
for the large volume limit. For example, both skyscraper sheaves of points
and $\mu$-stable vector bundles are stable;
the prediction that their phases differ by $\frac n2$ is reflected
by $\phi(\OO_x)(+\infty) = 1$ and
$\phi(\EE)(+\infty) = 1 - \frac n2$.

It may be worth mentioning that even for vector bundles and $\beta = 0$,
stability at the large volume limit does not coincide with Gieseker-stability.
Both stability conditions are refinements of slope stability, but 
they are different refinements.

If $X$ is a smooth Calabi-Yau variety, and if $2\beta = c_1(\LL)$ is the
first Chern class of a line bundle $\LL$, then the stability condition
is self-dual in the sense of Proposition \ref{prop:dual-stability}, 
with respect to $\LL^{-1}[n]$ as dualizing complex.

Now consider the case of a smooth projective surface.
Then $\AA^p$ is the category of two-term complexes
complexes $E$ with $H^{-1}(E)$ being torsion-free, and $H^0(E)$ being a
torsion sheaf. 

In the case of a K3 surface, the ample chamber is described
completely by \cite[Proposition 10.3]{Bridgeland:K3}; and for an
arbitrary smooth projective surface, the stability condition constructed in
\cite[section 2]{Aaron-Daniele} are also part of the ample chamber.
The following proposition gives a precise meaning to the catch phrase
``polynomial stability conditions at the large volume limit are limits
of Bridgeland stability conditions in the ample chamber'':
\begin{Prop-s} 				\label{prop:large-volume}
Let $S$ be a surface, $\beta \in A^1(X)_\R$ be a divisor class,
$\omega \in A^1(X)_\Q$ a \emph{rational} ample class, and let
$\rho, p$ be as above. Consider either of the following situations:
\begin{enumerate}[label={(\alph*)}]
\item \label{case:K3}
$S$ is a K3 surface; let $(Z_m, \PP_m)$ be the stability
condition constructed in \cite{Bridgeland:K3} from $\beta$ and
$\omega = n \cdot \omega_0$ (assuming $\omega^2 > 2$), and let
$(Z, \PP)$ be the polynomial stability condition constructed
from the data $\Omega$.
\item \label{case:surface}
$S$ is a smooth projective surface; let $(Z_m, \PP_m)$ be the
stability condition constructed in \cite{Aaron-Daniele} from
$\beta$ and $\omega = n \cdot \omega_0$, and let $(Z, \PP)$ be the 
polynomial stability condition constructed from
$\Omega' = (\omega, \rho, p, U=e^{-\beta})$.
\end{enumerate}
Then $E \in D^b(S)$ is
$(Z_m, \PP_m)$-stable for $m \gg 0$ if and only if it is $(Z, \PP)$-stable.
If $E \in D^b(S)$ is an arbitrary object, then the HN-filtration of
$E$ with respect to $(Z, \PP)$ is identical to the HN-filtration
with respect to $(Z_m, \PP_m)$ for $m \gg 0$.
\end{Prop-s}

In either case, the stability function is of the form
\begin{equation} 			\label{eq:surfaceZ}
 Z(E)(m) =  \ch_0(E) \omega^2 \cdot \frac{m^2}2
+ i \bigl(\omega \ch_1(E) - \ch_0(E) \beta \omega\bigr) m + c(E)
\end{equation}
for some real constant $c(E)$. Let
$\mu_\omega = \frac{\ch_1(E) \cdot \omega}{\ch_0(E)}$ be the slope
function for torsion-free sheaves on $S$ defined by $\omega$.

\begin{Lem-s}
Let $E \in D^b(S)$ be a $(Z, \PP)$-semistable object with
$0 \prec \phi(E) \preceq 1$.
Then $E$ satisfies one of the following conditions:
\begin{enumerate}[label={(\alph*)}]
\item \label{case:torsionsheaf}
$E$ is a $\mu_\omega$-semistable torsion sheaf.
\item \label{case:torsionfree}
$E$ is a torsion-free $\mu_\omega$-semistable sheaf
with $\mu_\omega(E) > \beta \cdot \omega$.
\item \label{case:torsionfreeshift}
$H^{-1}(E)$ is torsion-free $\mu_\omega$-semistable sheaf
of slope $\mu_\omega(H^{-1}(E)) \le \beta \cdot \omega$, 
$H^0(E)$ is zero-dimensional, and all other cohomology sheaves vanish.
\end{enumerate}
\end{Lem-s}
\begin{Prf}
Note that such an $E$ satisfies
$E \in \AA^p$ or $E \in \AA^p[-1]$, as $\AA^p = \PP((\frac 14, \frac 54])$.

If $E \in \AA^p[-1]$, then $H^1(E)$ is a torsion sheaf
by the definition of $\AA^p$. In fact, $H^1(E)$ has to vanish: otherwise
$\phi(H^1(E)) \preceq 1$, and because of
$\phi(E[1]) = \phi(E) + 1 \succ 1$ the
surjection $E[1] \onto H^1(E)$ would destabilize $E[1]$ in $\AA^p$. Hence
$E$ is a torsion-free sheaf.
Further, $E$ must be $\mu_\omega$-semistable: for
any surjection $E \onto B$ with $B$ torsion-free
and $\mu_\omega(E) > \mu_\omega(B)$, the surjection
that $E[1] \onto B[1]$ would destabilize $E[1]$ in $\AA^p$. Since
$\phi(E[1]) \succ 1$, we must have $\Im(Z(E)(m) < 0$ for $m \gg 0$;
this is equivalent to 
$\omega \ch_1(E[1]) - \ch_0(E[1]) \beta \omega < 0$ or
$\mu_\omega(E) > \beta \cdot \omega$. 

Similarly, one shows that if $E \in \AA^p$ and $H^{-1}(E)$ does not
vanish, then it is torsion-free and $\mu_\omega$-semistable of slope
$\mu_\omega (H^{-1}(E)) \le \beta \cdot \omega$. 
Also, $H^0(E)$ is of dimension zero:
otherwise $\phi(H^0(E))(+\infty) = \frac 12$, in contradiction to
$\phi(E)(+\infty) = 1$ and the surjection $E \onto H^0(E)$ in $\AA^p$.

Finally, if $E \in \AA^p$ and $H^{-1}(E)$ vanishes, then $E$
is a torsion sheaf, which is easily seen to be $\mu_\omega$-semistable.
\end{Prf}

\begin{figure}[htb]
\begin{center}
        \includegraphics{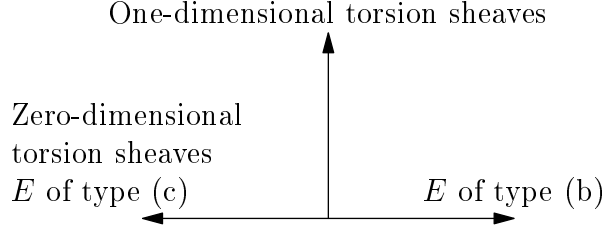}
\caption{Asymptotic directions of $Z(E)$ for $Z$-stable objects
$E \in \AA(\beta, \omega)$}
\label{fig:surface-stable}
\end{center} \end{figure}

\begin{Prf}[of proposition \ref{prop:large-volume}]
Let $\AA = \PP((0,1])$.
We first show that $\AA$ is identical to the heart
$\AA(\beta, \omega)$ defined in \cite[Lemma 6.1]{Bridgeland:K3},
respectively $\AA^\sharp_{(D, F)}$ defined in \cite[section 2]{Aaron-Daniele}.
Recall that $\AA(\beta, \omega)$ is characterized as the extension-closed
subcategory of $D^b(S)$ generated by
torsion sheaves, by $\mu_\omega$-semistable sheaves $F$ of
slope $\mu_\omega(F) > \beta \cdot \omega$, and by the shifts
$F[1]$ of $\mu_\omega$-semistable sheaves $F$ of slope 
$\mu_\omega(F) \le \beta \cdot \omega$.

Since $\AA(\beta, \omega)$ is extension-closed  and every $E$ in the
above list is an element of $\AA(\beta, \omega)$, it follows that
$\AA \subset \AA(\beta, \omega)$.  As both categories are
the heart of a bounded t-structure, they must be equal.

The first statement of the Proposition thus simplifies to
the claim that an object $E \in \AA$ is $Z$-stable
if and only if $E$ is $Z_m$-stable for $m \gg 0$.
By definition, we have $\phi(E) \succ \phi(F)$ if and only if
$\phi_m(E) = \phi(E)(m) > \phi_m(F)=\phi(F)(m)$
for $m \gg 0$; in particular,
if $E \in \AA$ is $Z$-unstable, then it will be $Z_m$-unstable for $m \gg 0$.

Conversely, assume that $E$ is $Z$-semistable.
In case \ref{case:torsionsheaf} of the lemma, $E$ is $Z_m$-stable for all $m$.
We now assume case \ref{case:torsionfreeshift}; case \ref{case:torsionfree}
can be dealt with similarly.
We need to show the following:
\emph{Given $E$, there is a constant $M$ such that whenever
$A \into E \onto B$ is a short exact sequence in $\AA$, 
then $\phi(E)(m) \le \phi(B)(m)$ for all $m \ge M$.}

If $B$ is a zero-dimensional
torsion sheaf, the claim is evidently satisfied. Otherwise
write $\EE:= H^{-1}(E)$,  $\BB := H^{-1}(B)$; let $\FF$ be the image of 
the induced map $\EE \to \BB$, and $\GG$ the cokernel.
The induced map $\GG \into H^0(A)$ has zero-dimensional cokernel;
hence $\beta \cdot \omega < \mu_\omega(H^0(A)) = \mu_\omega(\GG)$.
Since $\EE$ surjects
onto $\FF$, we have $\mu_\omega(\EE) \le \mu_\omega(\FF)$.
Combined with the definition of $\AA(\beta, \omega)$, we obtain
\begin{equation} 		\label{eq:bound-mu-B}
\mu_\omega(\EE) \le \mu_\omega(\BB) \le \omega \cdot \beta.
\end{equation}
Since $Z(B)(m)$ and $Z(E)(m)$ are in the semi-closed upper half plane
$\H$ for all $m$, the assertion is equivalent to
$\Im \bigl(\overline{Z(E)(m)} Z(B)(m)\bigr) \le 0$.
Using equation (\ref{eq:surfaceZ}) with
$\ch_0(E) = -\rk(\EE)$ and $\ch_1(E) = - \ch_1(\EE)$ etc.,
this can be simplified to:
\[
 \frac{\omega^2 m^2}2 \bigl(\mu_\omega(\BB) - \mu_\omega(\EE) \bigr)
 \ge \frac{c(B)}{\rk(\BB)}\bigl(\beta \omega - \mu_\omega(\EE) \bigr) \\
 - \frac{c(E)}{\rk(\EE)} \bigl(\beta \omega - \mu_\omega(\BB)\bigr)
\]
By inequality (\ref{eq:bound-mu-B}), all the expressions in parentheses
are non-negative. 
Since $E$ is $Z$-semistable, the inequality is satisfied for $m \gg 0$;
in particular, in the case $\mu_\omega(\BB) = \mu_\omega(\EE)$ it
holds for all $m$.
Excluding this case, the claim follows if we can bound
$\mu_\omega(\BB) - \mu_\omega(\EE)$ from below by a positive
constant and $\frac{c(B)}{\rk(\BB)}$ from above.

If $\GG = 0$, then the rank of $\BB$ is bounded above. By 
the rationality of $\omega$, the set of possible values of
$\omega \cdot \ch_1(\BB)$ is discrete,
giving a positive lower bound for $\mu_\omega(\BB) - \mu_\omega(\EE)$.
Otherwise, the lower bound follows from
$\mu_{\omega(\GG)} > \beta \cdot \omega$ and the upper bound on the rank of
$\FF$.

To prove the upper bound of $\frac{c(B)}{\rk(\BB)}$,
we restrict to the case \ref{case:surface} of the proposition.
Case \ref{case:K3}
can be proved similarly (and similarly to the proof of
\cite[Proposition 14.2]{Bridgeland:K3}); the argument is similar to
the proof of the existence of stability conditions in 
\cite[section 2]{Aaron-Daniele}.

It is sufficient to bound the number $\frac{c(B_j)}{\rk(\BB_j)}$
for every HN filtration quotient $B_j$ of $B$ with respect to $Z$,
and $\BB_j = H^{-1}(B_j)$. Then
$\BB_j$ is $\mu_\omega$-semistable, and its slope still satisfies the
inequality
\begin{equation} \label{eq:bound-mu-Bj}
\mu_\omega(\EE) \le \mu_\omega(\BB_j) \le \beta \cdot \omega.
\end{equation}

Using the Bogomolov-Gieseker inequality for
$\ch_2(\BB_j)$, we get:
\begin{align*}
c(B_j) 	& = -e^{-\beta} \cdot \ch(B_j) 
     	  = \ch_2(\BB_j) - \ch_2(H^0(B_j)) - \beta \cdot \ch_1(\BB_j)
		+ \rk(\BB_j) \frac{\beta^2}{2}  \\
	& \le \frac{\ch_1(\BB_j)^2}{2 \cdot \rk(\BB_j)}
		- \beta \cdot \ch_1(\BB_j) + \rk(\BB_j) \frac{\beta^2}{2}   \\
\frac{c(B_j)}{\rk(\BB_j)} & \le
	   \frac{1}{2} \left(\frac{\ch_1(\BB_j)}{\rk(\BB_j)} - \beta\right)^2
\end{align*}
Due to inequality (\ref{eq:bound-mu-Bj}) and the Hodge index theorem,
this number is bounded from above.

It remains to show the statement about the Harder-Narasimhan filtrations.
It is enough to show this for $E \in \AA$, as we already showed
$\PP((0,1]) = \AA = \PP_m((0,1])$. Let $A_1, A_2, \dots, A_n$ be the
Harder-Narasimhan filtration quotients of $E$ with respect to $\PP$.
Then the claim follows if $m$ is big enough such that every
every $A_j$ is $Z_m$-semistable, and such that
$\phi(A_1)(m) > \phi(A_2)(m) > \dots > \phi(A_n)(m)$.
\end{Prf}

\section{$X$ as the moduli space of stable point-like objects}

As a toy example of moduli problems in the derived category we will show that
in the smooth case, the moduli space of stable point-like objects is given by
$X$ itself. It shows that all our polynomial stability conditions are
``ample'' in the sense of the ample chamber in the introduction.

Given a polynomial stability condition $(Z, \PP)$ on $X$, a family of
stable objects over $S$ is an object $E \in D^b(X \times S)$ such that for
every closed point $s \in S$, the object $\L i_s^* E \in D^b(X)$ is 
$(Z, \PP)$-stable. Since $\Ext^{<0}(E, E) = 0$ for any stable object, it
is known that the moduli problem of stable objects is an abstract stack
(see \cite[Proposition 2.1.10]{Lieblich:mother-of-all} for a precise
statement and references). However, in general it is not known whether 
this stack is an algebraic Artin stack; see \cite{Toda:K3} for a proof in
a large class of examples.

Let $c$ be a class in the numerical $K$-group.
By some abuse of notation, we denote by $M_{c}(Z, \PP)$ the substack
of $(Z, \PP)$-stable objects such that $\L i_s^* E$ is an element of
$\AA^p$ and of class $c$.
\begin{Prop-s} \label{prop:X-as-modspace}
Assume that $X$ is a smooth projective variety over $\C$.
Let $(Z, \PP)$ be any of the polynomial stability conditions constructed
in Theorem \ref{mainthm} that has $p(0) = 0$.
The moduli stack $M_{[\OO_x]}(Z, \PP)$ of stable objects of the class
of a point is isomorphic to the trivial $\C^*$-gerbe $X/\C^*$ over $X$.
\end{Prop-s}
The assumption ensures that every skyscraper sheaf $\OO_x$ is 
an objects of $\AA^p$ (otherwise the same would be true after a shift,
and we might have to replace $[\OO_x]$ by $-[\OO_x]$ in the proposition).

\begin{Prf}
If $A \into \OO_x \onto B$ is a short exact sequence in $\AA^p$, then
the long exact cohomology sequence combined with lemma \ref{lem:dimsupport}
shows that $H^k(A) = 0 = H^k(B)$ for $k \neq 0$, and so $A \cong \OO_x$
or $B \cong \OO_x$. Hence every $\OO_x$ is stable.

Conversely, let $E \in \AA^p$ be any object with $[E] = [\OO_x]$. From lemma
\ref{lem:dimsupport} it follows that $H^k(E) = 0$ for $k \neq 0$, and
hence $E \cong \OO_x$ for some $x \in X$.

Hence the map $X \to M_{[\OO_x]}(Z, \PP)$ given by the structure sheaf
of the diagonal in
$X \times X$ is bijective on closed points. By the deformation theory
of complexes (see \cite[section 3]{Lieblich:mother-of-all} or
\cite{Inaba}) and $T_x X \cong \Ext^1(\OO_x, \OO_x)$, it induces an
isomorphism on tangent spaces. Since $X$ is smooth, the map is surjective.
\end{Prf}

\section{Wall-crossings: PT/DT-correspondence and one-dimensional
torsion sheaves}
\label{sect:PTDT}

In \cite{PT1}, Pandharipande and Thomas introduced new invariants of smooth
projective threefolds. They are obtained from moduli spaces of stable pairs
constructed by Le Potier in \cite{Potier:stable_pairs}; in their
context, a stable pair is a section $s \colon \OO_X \to \FF$ of a
pure one-dimensional sheaf $\FF$ that generically generates $\FF$.

In the Calabi-Yau case, the authors conjecture that 
the generating function of stable pairs invariants equals the reduced
generating function of Donaldson-Thomas invariants introduced in \cite{MNOP1}.
A heuristic justification of the conjecture was given in
\cite[section 3.3]{PT1} by interpreting the formula as a wall-crossing
formula under a change of Bridgeland stability conditions, assuming the
existence of certain stability conditions.

With proposition \ref{prop:PTDT}, we will show that this wall-crossing can
actually be achieved in a family of polynomial stability conditions, thus
making the heuristic justification one step more rigorous.

Further, in the subsequent article \cite{PT-BPS}, the authors give a new
geometric definition of BPS state counts. It relies on a relation between
invariants of stable pairs and invariants of one-dimensional torsion sheaves
(see \cite[Proposition 2.2]{PT-BPS}). In section \ref{sect:PT-BPS}, we show
that this relation can similarly be interpreted as a wall-crossing in our
family of polynomial stability condition; in fact the wall-crossing formula
is much simpler than in the case of the PT/DT-correspondence.

We refer to \cite{Toda:bir-BPS-states} for a similar use of a wall-crossing
to relate (differently defined) BPS state counts on birational Calabi-Yau
threefolds.

\subsection{PT/DT-correspondence}

Let $X$ be a smooth complex threefold. Fix an ample class $\omega \in A^1_\R$,
and let $p$ be
the perversity function $p(d) = -\lfloor \frac d2 \rfloor$.
Then the category of perverse coherent sheaves $\AA^p$ can be described
explicitly: a complex $E \in D^b(X)$ is an element of $\AA^p$ if
\begin{itemize*}
\item $H^i(E) = 0$ for $i \neq 0, -1$,
\item $H^0(E)$ is supported in dimension $\le 1$, and
\item $H^{-1}(E)$ has no torsion in dimension $\le 1$.
\end{itemize*}

\begin{figure}[htb]
\begin{center}
        \subfloat[DT-stability]{ \includegraphics{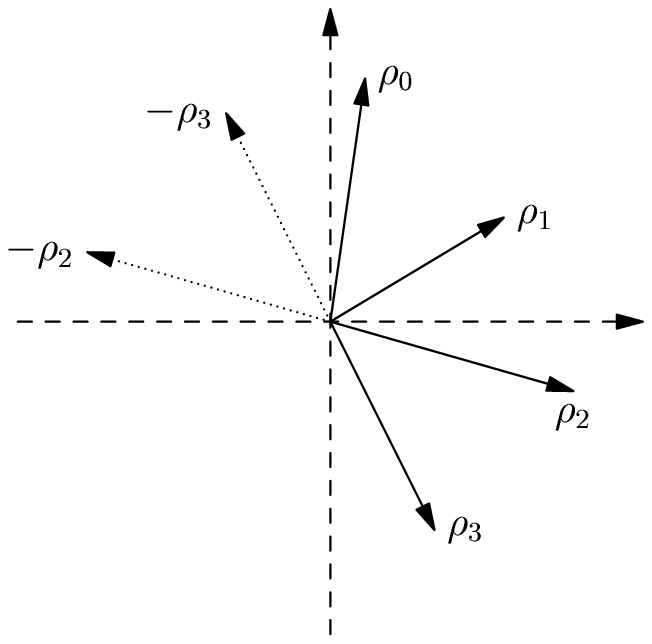} }
        \quad
        \subfloat[PT-stability]{ \includegraphics{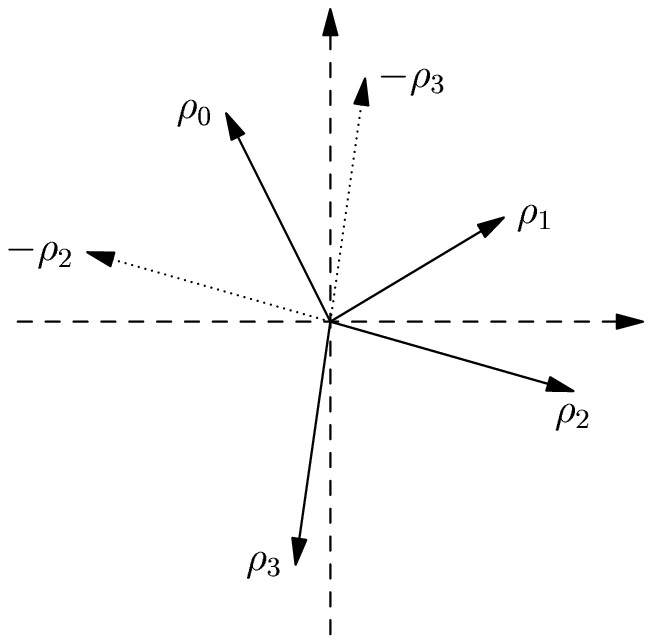} }
\caption{PT/DT wall-crossing}
\label{fig:rhos-PTDT} 
\end{center} \end{figure}

Consider stability vectors $\rho$ such that $p$ is an associated
perversity function, i.e. 
$\rho_0, \rho_1 \in \H$ and $\rho_2, \rho_3 \in -\H$.
Let $U$ be arbitrary, and consider the polynomial stability functions given
by \[ Z(E)(m) = \sum_{d=0}^3 \rho_d m^d \omega^d \cdot \ch(E)\cdot U. \]
We further assume $\phi(-\rho_3) > \phi(\rho_1)$.
We call it a \emph{DT-stability function} if
$\phi(-\rho_3) > \phi(\rho_0)$ and a \emph{PT-stability function}
if $\phi(\rho_0) > \phi(-\rho_3)$, see figure \ref{fig:rhos-PTDT}.

If $\omega$ is the class of an ample line bundle $\LL$, $U = \td X$ and
$h(E)(m) = \sum_{d=0}^3 a_d m^d$ is the Hilbert polynomial of $E$ with
respect to $\LL$, then the central charges can also be written as
the complexified Hilbert polynomial
$Z(E)(m) = \sum_{d=0}^3 d!\rho_d a_d m^d$. 

Fix numerical invariants $\beta \in A_1^{\num}$ and $n \in A_0 \cong \Z$.
We consider the moduli problem
$M_{(-1,0,\beta,n)}^0(Z, \PP)$ of $Z$-stable objects in 
$\AA^p$ with trivialized determinant and numerical invariants in $A_*^{\num}$
given by $\ch(E) = (-1, 0, \beta, n)$.

\begin{Prop} 			\label{prop:PTDT}
Let $S$ be of finite type over $\C$, and $I \in D^b(X \times S)$
be an object with $\ch(I_s) = (-1, 0, \beta, n)$ for every closed point
$s \in S$, and with trivialized determinant.

If $Z$ is a DT-stability function, then $I$
is a $Z$-stable family of objects in $\AA^p$
if and only if it is quasi-isomorphic to the shift $\JJ[1]$ of a flat
family of ideal sheaves of one-dimensional subschemes.

If $Z$ is a PT-stability function and $\beta \neq 0$,
then $I \in D^b(X \times S)$ is a $Z$-stable family of objects
in $\AA^p$ if and only if it is
quasi-isomorphic to the complex $\OO_{X \times S} \to \FF$
(with $\FF$ in degree zero) of a family of stable pairs as defined
in \cite{Potier:stable_pairs, PT1}.
\end{Prop}
Thus in both cases we get an isomorphism
of moduli spaces of ideal sheaves/stable pairs with the moduli space
$M_{(-1,0,\beta,n)}^0(Z, \PP)$ of $Z$-stable objects with trivialized
determinant.

If $\beta = 0$ and $Z$ is a PT-stability function, then the only semistable
object is $\OO_X[1]$ of class $(-1, 0, 0, 0)$. This does not agree with the
definition of stable pairs, but does give the correct generating
function, so that the conjectured wall-crossing formula of \cite{PT1} holds
for all $\beta$.

\begin{Prf}
Let $Z$ be a DT-stability function, and let assume that $I$ is
a family of stable objects. If for any closed point $s \in S$,
we would have both
$H^{-1}(I_s) \neq 0$ and $H^0(I_s) \neq 0$, then the short exact sequence
\[ H^{-1}(I_s)[1] \to I_s \to H^0(I_s) \]
would destabilize $I_s$: for large $m$, the phase of
$Z(H^{-1}(I_s)[1])(m)$ is approaching $\phi(-\rho_2)$ or $\phi(-\rho_3)$
(depending on the dimension of support of $I_s$);
while the phase of $Z(H^0(I_s))(m)$ is approaching 
$\phi(\rho_1)$ or $\phi(\rho_0)$. Hence $I$ is the shift of a flat family
$\JJ$ of sheaves of rank one.\footnote{Here, and again later in the proof of
the PT-case, we are using the following standard fact (cf. \cite[Lemma
3.31]{Huybrechts:FM}): If $I$ is a complex on $X \times S$ such that for every
closed point $s \in S$, the derived pull-back $I_s$ is a sheaf, then $I$ is a
sheaf, flat over $S$.}
To be both stable and an element of $\AA^p$, it has to be torsion-free.
Its double dual is locally free by \cite[Lemma 6.13]{Kollar:projective_moduli}.
Since $\JJ$ has trivialized determinant, the double dual $\JJ^{**}$ is
the structure sheaf $\OO_{X \times S}$; the natural inclusion
$\JJ \into \JJ^{**}$ exhibits $\JJ$ as a flat family of ideal sheaves.
Conversely, any such flat family of ideal sheaves gives a family of
stable objects in $\AA^p$.

(In fact, the DT-stability conditions are obtained from the stability
conditions of section \ref{sect:Simpson} corresponding to Simpson stability by
a rotation of the complex plane and accordingly tilting the heart of the
t-structure. Hence the stable objects are exactly the shifts of Simpson-stable
sheaves; their moduli space is well-known to be isomorphic to the Hilbert
scheme.)

Now let $Z$ be a PT-stability function.
We have to show that $I_s$ is stable for all $s$ if and only if
$I$ is quasi-isomorphic to a family of stable pairs: a complex
$\OO_{X \times S} \to \FF$ such that 
\begin{enumerate}
\item $\FF$ is flat over $S$, and
\item $\OO_X \to \FF_s$ is a stable pair for all $s \in S$.
\end{enumerate}
Given such a family of stable pairs, the associated complex
is a family of objects in $\AA^p$ with trivialized determinant.

First assume that $I$ is $Z$-stable. 
By the same argument as in the
DT-case, $H^0(I_s)$ must be zero-dimensional, and
$H^{-1}(I_s)[1]$ torsion-free of rank one with trivialized determinant.

It follows that $Q:= H^0(I)$ is zero-dimensional over $S$, and that
$H^{-1}(I)$ is a torsion-free rank one sheaf with trivialized determinant.
Let $U \subset X \times S$ be the complement of the support of $Q$, and
let $I_U := I|_{U}[-1]$ be the restriction of $I[-1]$ to $U$. Then the
derived pull-back of $I_U$ to every fiber over $s \in S$ is a
sheaf; so $I_U \cong H^{-1}(I)|_U$ is itself a sheaf,
flat over $S$.
Hence $H^{-1}(I)$ is flat over $S$ outside a set of codimension 3.

By the same arguments as in the proof of \cite[Theorem 2.7]{PT1} it follows
that $H^{-1}(I)$ is a family of ideal sheaves $\JJ_Z$ of one-dimensional
subschemes of $X$.  The complex $I$ is the cone of a map $Q \to \JJ_Z[2]$.
Since $Q$ is zero-dimensional over $S$, we have
$\Ext^1(Q, \OO_{X \times S}) = 0 = \Ext^2(Q, \OO_{X \times S})$; combined
with the short exact sequence
$\JJ_Z \to \OO_{X \times S} \to \OO_Z$, we get a unique factorization
$Q \to \OO_Z[1] \to \JJ_Z[2]$. Using the octahedral axiom associated to
this composition, we see that $I$ is the cone of a map
$\OO_{X \times S} \to \FF$,
where $\FF$ (in degree zero) is the extension of $\OO_Z$ and $Q$ given as the
cone of the map $Q \to \OO_Z[1]$ above.

It remains to prove that $I_s$ is $Z$-stable if and only if 
$\OO_X \to \FF_s$ is a stable pair. Assume that $I_s$ is
$Z$-stable, and
note that $\phi(I_s)(+\infty) = \phi(-\rho_3)$.

Since $\beta \neq 0$, the sheaf $\FF_s$ is one-dimensional.
It cannot have a zero-dimensional subsheaf
$Q \into \FF_s$, as this would induce an inclusion
$Q \into I_s$ in $\AA^p$, destabilizing $I_s$ due to
$\phi(Q) = \phi(\rho_0) > \phi(-\rho_3)$.
Thus $\FF_s$ is purely one-dimensional, and $\OO_X \to \FF_s$ is stable
by \cite[Lemma 1.3]{PT1}.

Conversely, assume that $I_s$ is a stable pair. Consider any destabilizing
short exact sequence $A \into I_s \onto B$ in $\AA^p$ with
$\phi(A) \succ \phi(I_s) \succ \phi(B)$, and its long exact
cohomology sequence 
\[ H^{-1}(A) \into H^{-1}(I_s) \to H^{-1}(B) \to
   H^0(A) \to H^0(I_s) \onto H^0(B). \]
If $H^{-1}(A) \into H^{-1}(I_s)$ is a proper inclusion, then $H^{-1}(B)$ is
supported in dimension 2, and we get the contradiction $\phi(B)(+\infty) =
\phi(-\rho_2) > \phi(I_s)(+\infty)$. So either $H^{-1}(A) = H^{-1}(I_s)$ or
$H^{-1}(A) = 0$. In the former case, $H^{-1}(B) = 0$; since $B = H^0(B)$ is
supported in dimension zero, we get the contradiction
$\phi(B) = \phi(\rho_0) > \phi(I_s)(+\infty)$.  In the latter case,
$A = H^0(A)$ must be zero-dimensional to destabilize $I_s$; by the purity of
$\FF_s$, this implies $\Hom(A, \FF_s) = 0$.  Together with 
$\Ext^1(A, \OO_X) = 0$ and the exact triangle $\FF_s \to I_s \to \OO_X[1]$,
this shows the vanishing of $\Hom(A, I_s)$.
\end{Prf}

The reason to expect a wall-crossing formula in a situation as above is the
following: Denote by $Z_{\PT}$ a PT-stability function, and by $Z_{\DT}$ a
DT-stability function. If $E$ is $Z_{\PT}$-semistable but $Z_{\DT}$-unstable,
then we can write $E$ as an extension of $Z_{\DT}$-semistable objects (by the
existence of Harder-Narasimhan filtrations); and conversely for
$Z_{\DT}$-semistable but $Z_{\PT}$-unstable objects.  Hence one can expect an
expression for the difference between the counting invariants of $Z_{\DT}$-
respectively $Z_{\PT}$-semistable objects in terms of lower degree counting
invariants. This observation (due to D. Joyce, cf. \cite{Joyce4}) can be
made more concrete and precise in the situation considered below.

\subsection{Stable pairs and one-dimensional torsion sheaves}
\label{sect:PT-BPS}

Let $X$ be a Calabi-Yau threefold, and $\beta, n$ as before. By a counting
invariant we will always denote the signed weighted Euler characteristic
(in the sense of \cite{Behrend:DT-microlocal}) of a moduli space of
stable objects of some fixed numerical class, and with trivialized
determinant.

In the very recent preprint \cite{PT-BPS}, the authors give a new
geometric description of BPS state counts for irreducible curve classes
on $X$. They use the counting invariants 
$N_{n, \beta}$ of the moduli spaces $\MM_n(X, \beta)$ of stable
one-dimensional torsion sheaves of class $(0, 0, \beta, n)$.
At the core of their argument is the following relation: if
$\beta$ is an irreducible effective class and
$P_n(X, \beta)$ denotes the counting invariant of stable pairs of
class $(-1, 0, \beta, n)$, they prove that
\begin{equation}			\label{eq:BPS-relation}
P_n(X, \beta) - P_{-n}(X, \beta) = (-1)^{n-1}n N_{n, \beta}.
\end{equation}

To make the subsequent discussion more specific, 
we fix $\rho_0 \in \R_{>0}\cdot (-1)$, $\rho_1 \in \R_{>0}\cdot i$, 
$\rho_2 \in \R_{>0}$.  We keep
$\omega, p$, and in particular continue
to work with the same category of perverse coherent sheaves $\AA^p$.
Assume that $P(\overline U) = U$, e.g. $U \in A^{\text{even}}(X)_\R$.
For $a > 0$ write
$Z_a$ for the polynomial stability function on $\AA^p$ obtained from 
$\rho_3 = - b \cdot i + a$ (for some $b > 0$), and similarly
$Z_{-a}$ for $\rho_3 = - b \cdot i - a$; see also the figure.

\begin{floatingfigure}[htb]{160pt}
\begin{center}
        \includegraphics{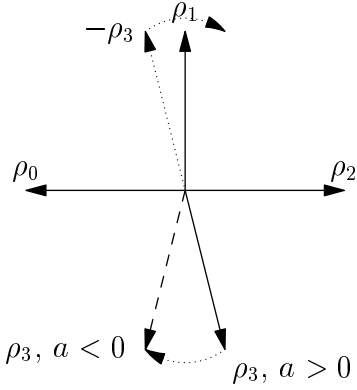}
\caption{Wall-crossing between PT invariants and BPS state counts}
\label{fig:PT-BPS}
\end{center}
\end{floatingfigure}

Then $Z_a$ is a ``PT-stability function'' (in the terminology of the
previous section), hence the stable objects of
class $(-1, 0, \beta, n)$ are the stable pairs
$I \cong \OO_X \to \FF$. If we cross the wall $a = 0$ (the large volume
limit), then short exact sequence
$\FF \to I \to \OO_X[1]$ destabilizes $I$ for $a < 0$.

Let $\D$ be the dualizing functor $E \mapsto \RlHom(E, \OO_X[2])$.
Then the polynomial stability condition obtained from $Z_a$ is dual
to that of $Z_{-a}$; this can be seem from proposition \ref{prop:dual-stability}
and the fact that in our case $\AA^p$ is a tilt of $\AA^{p^*}$,
compatible with the stability condition.

It follows that if $a < 0$, then the stable objects of same class are
the derived duals of stable pairs of class $(-1, 0, \beta, -n)$; their
counting invariant is thus given by $P_{-n}(X, \beta)$. 

If we additionally assume that $\beta$ is irreducible, then $\FF$ is
stable for both $Z_a$ and $Z_{-a}$, and
the short
exact sequence $\FF \to I \to \OO_X[1]$ is the HN filtration 
with respect to $Z_{-a}$ of a stable pair $I$. Conversely, the dual 
short exact sequence $\OO_X[1] = \D(\OO_X[1]) \to \D(I) \to \D(F)$
will be the HN filtration of $\D(I)$ with respect to $Z_a$ (where
$\D(F)$ is a stable sheaf of class $(0, 0, \beta, n)$).
Hence the wall-crossing formula
can be written schematically as
\begin{align*} P_n(X, \beta) - P_{-n}(X, \beta) 
= & \sharp\bigl(\text{Extensions of $\OO_X[1]$ with $\FF$}\bigr) \\
 & -\sharp\bigl(\text{Extensions of $\FF'$ with $\OO_X[1]$}\bigr),
\end{align*}
where $\FF, \FF'$ can be any stable sheaf of class $(0, 0, \beta, n)$.

If the dimensions of $\Ext^1(\OO_X[1], \FF) = H^0(\FF)$ and $\Ext^1(\FF',
\OO_X[1]) = H^1(\FF')^*$ were constant, then the moduli spaces of extensions
would be projective bundles over $\MM_n(X, \beta)$; in this case, formula
(\ref{eq:BPS-relation}) would follow immediately.
Without this simplifying assumption, one can still hope 
to prove formulas such as (\ref{eq:BPS-relation}) using a stratification
of the moduli spaces and the formalism of \cite{Behrend:DT-microlocal}.
In fact, the proof in \cite{PT-BPS} exactly follows this general principle,
the key ingredient being a control of the constructible
functions of \cite{Behrend:DT-microlocal} by \cite[Theorem 3]{PT-BPS}.

\section{Existence of Harder-Narasimhan filtrations}
\label{sect:HN-proof}

In this section we will prove that the category of perverse coherent
sheaves has the Harder-Narasimhan property for the polynomial stability
function $Z$ defined in Theorem \ref{mainthm}. The proof is complicated
by the fact that $\AA^p$ is in general neither $Z$-Artinian nor
$Z$-Noetherian.

\subsection{Perverse coherent sheaves and tilting}

Essential for the proof is a more detailed understanding of the category
of perverse coherent sheaves, more precisely the existence of certain
torsion pairs in that category.  We recall briefly the notion of a torsion
pair and a tilt of a t-structure:

\begin{Def} 
A torsion pair in an abelian category $\AA$ is a pair of full subcategories
$\TT, \FF$ such that 
\begin{enumerate}[label={(\alph*)}]
\item $\Hom(T, F) = 0$ for all $T \in \TT$ and $F \in \FF$.
\item For every $E \in \AA$ there is a short exact sequence
$T \into E \onto F$ in $\AA$ with $T \in \TT$ and $F \in \FF$.
\end{enumerate}
\end{Def}
If $(\TT, \FF)$ satisfy both conditions, then $T, F$ are uniquely
determined by $E$. They depend functorially on $E$, and the functors
$E \mapsto T$ and $E \mapsto F$ are left-exact and right-exact, respectively.

Now assume $\AA$ is the heart of a bounded t-structure in a triangulated
category $\DD$, with associated
cohomology functors $H_\AA^i \colon \DD \to \AA$. Given a torsion pair $\TT,
\FF$ in $\AA$, then the following defines the heart $\AA^\sharp$ of a
related t-structure (called the tilt of $\AA$; see \cite{Happel-al:tilting}): 
An object $A$ is in $\AA^\sharp$ if 
\[
H^0_\AA(A) \in \TT, \quad
H^{-1}_\AA(A)  \in \FF, \quad \text{and} \quad
H^{i}_\AA(A) = 0 \quad \text{if $i \neq 0, -1$.}
\]
The new heart $\AA^\sharp$ evidently satisfies
$\AA^\sharp \subset \langle \AA, \AA[1] \rangle$, and on the other hand
every heart of a bounded t-structure with this property is obtained as a 
tilt. This is shown by the following lemma, which is a slight reformulation
of a lemma in \cite{Polishchuk:families-of-t-structures}:
\begin{Lem}
Let $\AA, \AA^\sharp$ be the hearts of bounded t-structures in
a triangulated category $\DD$. If they satisfy
$\AA^\sharp \subset \langle \AA, \AA[1] \rangle$ (or, equivalently,
$\AA \subset \langle \AA^\sharp, \AA^\sharp[-1] \rangle$), then
\[
 \TT := \AA \cap \AA^\sharp, \quad
 \FF := \AA \cap \AA^\sharp[-1]
\]
defines a torsion pair in $\AA$, the heart $\AA^\sharp$ is obtained
from $\AA$ by tilting at this torsion pair, and $\FF[1], \TT$ is a
torsion pair in $\AA^\sharp$.
\end{Lem}
\begin{Prf}
If $(\DD^{\ge 0}, \DD^{\le 0})$ and 
$(\DD^{\sharp,\ge 0}, \DD^{\sharp,\le 0})$ are the two t-structures,
either assumption is equivalent to either of the following equivalent
assumptions:
\[
\DD^{\ge 0} \subset \DD^{\sharp, \ge 0} \subset \DD^{\ge -1}
	\quad \text{or} \quad
\DD^{\le 0} \supset \DD^{\sharp, \le 0} \supset \DD^{\le -1} 
\]
This is the assumption of
\cite[Lemma 1.1.2]{Polishchuk:families-of-t-structures}.
\end{Prf}

Now consider a perversity function $p$ and any $k \in \Z$ such that
$k = - p(d)$ for some $0 \le d \le n$. Consider the function
$p^k \colon \{0, \dots, n\}$ defined by
\[ p^k(d) = \begin{cases} p(d) \quad \text{if $p(d) \ge -k$} \\
			  p(d)+1 \quad \text{if $p(d) < -k$} \end{cases}
\]
Then $p^k(d)$ is a perversity function, and the hearts of
perverse coherent sheaves $\AA^{p^k}, \AA^p$ satisfy the assumptions of
the lemma. Hence 
\[ \FF_k = \AA^p \cap \AA^{p^k}, \quad \TT_k = \AA^p \cap \AA^{p^k}[1] \]
defines a torsion pair in $\AA^p$.

From the definition of the t-structures in Theorem
\ref{thm:tstruct}, and from lemma
\ref{lem:dimsupport}, it can easily be seen that the torsion pairs
can be described as below:
\begin{Prop}				\label{prop:TkFk}
Let $p$ be a perversity function and $k \in \Z$ such that
$p(0) \ge -k > p(n)$.
There is a torsion pair $(\TT_k, \FF_k)$ in $\AA^p$ defined
as follows:
\begin{flalign*}
\FF_k &= \stv{E \in \AA^p}{H^{-k'}(E) = 0 \quad \text{for $k' > k$}} \\
\TT_k &= \stv{E \in \AA^p}{H^{-k'}(E) \in \AA^{p, \le k'-1}
	\quad \text{for $k' \le k$}}
\end{flalign*}
The subcategory $\FF_k$ is closed under subobjects and quotients.
\end{Prop}
The only thing left to prove is the statement about $\FF_k$.
It is always the case for a torsion pair that $\FF$ is closed under subobjects
and $\TT$ under quotients. That $\FF_k$ is additionally closed under
quotients follows easily from the long exact cohomology sequence.

We denote by $\tau^k_\TT \colon \AA^p \to \TT_k$ and
$\tau^k_\FF \colon \AA^p \to \FF_k$ the associated functors; then for any
short exact sequence $A \into E \onto B$ in $\AA^p$ there is a (not very)
long exact sequence
\begin{equation}			\label{eq:Tk-Fk-sequence}
 \tau^k_\TT A \into \tau^k_\TT E \to \tau^k_\TT B \to 
   \tau^k_\FF A \to \tau^k_\FF E \onto \tau^k_\FF B.
\end{equation}

Note that $\tau^k_\FF$ will in general not coincide with
the truncation functor $\tau_{\ge -k}$ of the standard t-structure;
in fact, given $E \in \AA^p$ there is no reason why
$\tau_{\ge -k} (E)$ should also be an object of $\AA^p$.

\subsection{Dual stability condition}

The proof in the following section is substantially simplified by the use of
the dual stability condition constructed in Proposition
\ref{prop:dual-stability}. To use it, we need a partial proof of the duality
here. 

It is constructed from the dual t-structure. Let $\omega_X, D, \D$ be
as in section \ref{sect:dual-stability1}.
\begin{Prop}[\cite{Bezrukavnikov:perverse}]
Let $p$ be a perversity function, and $\barp$ the dual perversity
function (cf. definition \ref{def:perversity}); let 
$p^* = \barp + D - n$ be the dual perversity normalized according
to the choice of $\omega_X$. Define
$D^{p^*, \ge 0}, D^{p^*, \le 0} \subset \D(D^b(X))$
by the analogues of equations (\ref{eq:le0})
and (\ref{eq:ge0}), respectively.
Then the t-structures associated to $p, p^*$ are dual to each other
with respect to $\D$:
\[
 \D \left( D^{p, \le 0} \right) = D^{p^*, \ge 0} \quad \text{and} \quad
 \D \left( D^{p, \ge 0} \right) = D^{p^*, \le 0}
\]
\end{Prop}

By some abuse of notation, we will write $\AA^{p^*}$ for the intersection
$D^{p^*, \ge 0} \cap D^{p^*, \le 0} \subset \D(D^b(X))$.

\begin{Lem}
Given $\Omega$ and $\Omega^*$ as in Proposition \ref{prop:dual-stability},
$Z_{\Omega^*}$ is a polynomial stability function for the 
category of perverse sheaves $\AA^{p^*}$ of the dual perversity.
If $\phi, \phi^*$ are the polynomial phase functions of
$\AA^p, Z_\Omega$ and $\AA^{p^*}, Z_{\Omega^*}$, respectively,
then 
\begin{equation}			\label{eq:dualslopes}
 \phi(E_1) \prec \phi(E_2) \Leftrightarrow
\phi(\D(E_2)) \prec \phi(\D(E_1)).\end{equation}
An object $E \in \AA^p$ is $Z_{\Omega}$-stable if and only if
$\D(E) \in \AA^{p^*}$ is $Z_{\Omega^*}$-stable.
\end{Lem}

\begin{Prf}
Since $\ch(\D(E)) = P(\ch(E)) \cdot \ch(\omega_X)$, we have
\begin{align*}
Z_{\Omega^*}(\D(E)))(m)
&= \int_X \sum_{d=0}^n (-1)^{d+D}\overline{\rho_d} \omega^d m^d \cdot
P(\ch(E))\ch(\omega_X) \cdot (-1)^D \ch(\omega_X)^{-1} P(\overline{U})		\\
& = \int_X \sum_{d=0}^n m^d P(\overline{\rho_d} \omega^d)
P(\ch(E)) = (-1)^n \overline{Z(E)(m)}.
\end{align*}
This shows that $Z_{\Omega^*}$ is a polynomial stability function,
as $Z_{\Omega^*}(\D(E)(m))$ is in the interior of
$(-1)^{n+1} e^{-i\epsilon} \cdot \H$ whenever
$Z_{\Omega}(E(m))$ is in the interior of $e^{i\epsilon} \cdot\H$;
it also shows the equivalence (\ref{eq:dualslopes}).

Since $\D$ turns inclusions $E_1 \into E_2$ in $\AA^p$ into
quotients $\D(E_2) \onto \D(E_1)$ in $\AA^{p^*}$, and vice versa, this
also implies the claim about stable objects.
\end{Prf}

The lemma yields part (a) and (b) of Proposition \ref{prop:dual-stability}.

\subsection{Induction proof}

\begin{Lem}			\label{lem:quotcat}
Consider the quotient category
$\AA^{p, =k} = \AA^{p, \le k}/\AA^{p, \le k-1} \cong
\FF_k/\FF_{k-1}$
and let $Z' \colon K(\AA^{p, =k}) \to \C[m]$ be the restricted stability
function defined by 
\[ Z'(E)(m) = \int_X \sum_{\substack{d \in \{0, \dots, n\} \\ p(d) = -k}}
			\rho_d \omega^d m^d \cdot \ch(E) \cdot U
\]
Then $\AA^{p, =k}$ is Noetherian and strongly $Z'$-Artinian.
\end{Lem}
Here ``strongly $Z'$-Artinian'' says that there is no sequence of
inclusions
\[ \dots \into E_{j+1} \into E_j \into \dots \into E_2 \into E_1 \]
as in Proposition \ref{prop:HNproperty} with the weaker assumption
$\phi(E_{j+1}) \succeq \phi(E_j)$ for all $j$.

\begin{Prf}
For both statements, the proof is almost identical to the proof of the
same statement for $\AA = \Coh X$ and Simpson stability. We
will prove that the category is strongly $Z'$-Artinian.

Consider an infinite sequence of inclusions as above.
Since the dimension of the
support of $E_j$ is decreasing, we may assume it is constant, equal to
$d$. Similarly, we may assume that the lengths of $E_j$ at the generic
points of the (finitely many)
$d$-dimensional
components of its support are constant. In particular
the leading term of $Z'(E_j)(m)$ given by
$\omega^d \cdot \ch_{n-d}(E_j) \rho_d \cdot m^d$ is constant. 
The quotient $B_j = E_j/E_{j+1}$ is supported in strictly smaller dimension
$d' < d$.
Hence the leading term of $Z'(B_j)(m)$ is a positive linear multiple of
$\rho_{d'} m^{d'}$.
This implies
$\phi(E)(+\infty) = \phi(\rho_d) < \phi(\rho_{d'}) = \phi(B)(+\infty)$,
since $p(d') = p(d)$ and $p$ is a perversity function associated
to $\rho$.
Thus $\phi(E_j) \prec \phi(B_j)$, in contradiction to
$\phi_{E_{j+1}} \succeq \phi_{E_j}$ and the see-saw property.
\end{Prf}

We now come to the main proof.  As mentioned before, we can't apply
Proposition \ref{prop:HNproperty}.  Nevertheless, our proof follows
Bridgeland's proof of the corresponding statement \cite[Proposition
5.3]{Bridgeland:Stab} quite closely:
\begin{description*}
\item[Step 1] Every non-semistable $E \in \AA^p$ has a semistable
subobject $A \into E$
such that $\phi(A) \succ \phi(E)$, and a semistable quotient
$E \onto B$ with $\phi(E) \succ \phi(B)$.
\item[Step 2]
Every object $E$ has a maximal destabilizing quotient (mdq) $E \onto B$.
\item[Step 3]
Let $E_{j+1} \into E_j \into \dots \into E_1$ be the sequence
of inclusions in $\AA^p$
determined by $B_j$ being the mdq of $E_j$, and $E_{j+1}$ being the
kernel of the surjection $E_j \onto B_j$. Then this sequence terminates.
\end{description*}

A mdq is a quotient $E \onto B$ such that
for every other quotient $E \onto B'$, we have $\phi(B') \succeq \phi(B)$,
and such that equality holds if and only if the quotient factors via
$E \onto B \onto B'$. The proof of [ibid.] shows that the existence of
Harder-Narasimhan filtrations is equivalent to the existence of a mdq for
every object, and the termination of the sequence defined in step 3.

\subsubsection*{Step 1} 
Define a sequence of inclusions as follows: \emph{If $E_j$ is not semistable,
then among all subobjects $A \into E_j$ with $\phi(A) \succ \phi(E_j)$,
let $E_{j+1}$ be one such that the dimension of support $d(B_j)$ 
of $B_j$ is maximal, where
$B_j$ is the cokernel of $E_{j+1} \into E_j$.} It suffices
to prove that this sequence terminates.

By the definition of $E_{j+1}$, the sequence $d(B_j)$ of dimension
of support is monotone decreasing. By induction,
we just need to show that any such sequence
with $d(B_j) = d$ for all $j$ must terminate.

Let $k = -p(d)$, and consider the functors $\tau^k_\TT, \tau^k_\FF$ of
Proposition \ref{prop:TkFk}. Since $B_j \in \FF_k$, we have
$\tau^k_\TT (B_j) = 0$. By the exact sequence
(\ref{eq:Tk-Fk-sequence}), this shows that
$\tau^k_\TT (E_{j+1}) = \tau^k_TT (E_j)$ and that
\[ 0 \to \tau^k_\FF (E_{j+1}) \to \tau^k_\FF (E_j) \to B_j \to 0 \]
is exact. Taking cohomology, we get an induced short exact sequence
\[ 0 \to H^{-k}\left(\tau^k_\FF (E_{j+1})\right)
     \to H^{-k}\left(\tau^k_\FF (E_j)\right) \to H^{-k}\left(B_j\right) \to 0 \]
in $\AA^{p, =k}$.
From the lemma it follows that there must be a $j_0$ with 
\[ \phi\left(\tau^k_\FF (E_{j_0+1})\right) \prec \phi\left(\tau^k_\FF
(E_{j_0})\right) \prec \phi\left(B_{j_0}\right). \]
By the see-saw property, $\tau_\TT^k E_{j_0}$ is another subobject
of $E_{j_0}$ with $\phi(\tau_\TT^k E_{j_0}) \succ \phi(E_{j_0})$.
By the definition of $E_{j_0+1}$, this implies
$d(\tau^k_\FF (E_j)) = d(B_j)$ for $j= j_0$, and thus also for all
$j > j_0$, which is impossible.

This shows that every object has a semistable subobject as desired. By
applying the same arguments to the dual perversity and dual stability
function, this also shows that every object has a semistable quotient as
claimed.

\subsubsection*{Step 2}
We will prove steps 2 and 3 in a 2-step induction: To prove step 2
for an object supported in dimension $d$, we assume that steps 2 and 3
have been proven for objects supported in dimension at most $d-1$. To prove
step 3, we will assume that step 2 has been proven in dimension
$d$, and that step 3 has been proven in dimension $d-1$. The reason
this induction works well is that the subcategory of $\AA^p$ of objects
supported in dimension at most $d$ is closed under subquotients.

To prove step 2, we will instead show the dual statement:
Every object has a \emph{minimal
destabilizing subobject} (mds), i.e. a subobject $A \into E$ such that for 
every $A' \into E$ we have $\phi(A) \preceq \phi(A')$, with equality 
if and only if there is a factorization $A' \into A \into E$. 

Let $E_1 \in \AA^p$ be supported in dimension $d$, and let
$k = -p(d)$.  Define the sequence of objects $E_j$ as follows:
\begin{enumerate}
\item If $E_j$ is semistable, stop.
\item If there is a semistable
quotient $E_j \onto B_j$ with $\phi(E_j) \succ \phi(B_j)$ and
$H^{-k}(B_j) \neq 0$, then let $E_{j+1}$ be its kernel.
\item 
Otherwise,
let $B_j$ be the maximal destabilizing
quotient of $\tau^{k-1}_\FF E_j$, which exists by induction; and
$E_{j+1}$ be the kernel of the composition
$E_j \onto \tau^{k-1}_\FF E_j \onto B_j$.
\end{enumerate}
If neither case (1) nor (2) applies, there must be a
a semistable quotient $E \onto B$ with $\phi(E) \succ \phi(B)$ 
and $H^{-k}(B) = 0$. Then the quotient must factor as
$E \onto \tau^{k-1}_\FF E_j \onto B$. Then the mdq $B_j$ of
$\tau^{k-1}_\FF E_j$ satisfies $\phi(B) \succ \phi(B_j)$ by definition.

Hence both in case (2) and (3), we have a short exact sequence
$E_{j+1} \into E_j \onto B_j$ with $B_j$ semistable and
$\phi(E_j) \succ \phi(B_j)$.  By the arguments dual to those
given by Bridgeland, a mds of $E_{j+1}$ is also be a mds of $E_j$,
and if $E_j$ is semistable it is its own mds. So we just need to prove
that the above algorithm terminates.

By the lemma, case (2) will only happen a finite number of times.
However, in case (3) we get a short exact sequence
\[ \tau^{k-1}_\FF E_{j+1} \into \tau^{k-1}_\FF E_j \onto B_j, \]
where $B_j$ is the mdq of $\tau^{k-1}_\FF E_j$. By the induction
assumption about step 3, this sequence must terminate as well.

Finally, note that if $E$ is supported in dimension $d$, then so is
$\D(E)$. Again we can use the same arguments in the dual setting and prove the
existence of an mdq for objects supported in dimension $d$ as well.

\subsubsection*{Step 3}
Let $k = -p(\dim E_1)$. Again, by lemma \ref{lem:quotcat}, the sequence
of inclusions $H^{-k}(E_{j+1}) \into H^{-k}(E_j)$ will become an
isomorphism in the quotient category $\AA^{p, =k}$ after a finite
number of steps. Then $H^{-k}(B_j)$ is in
$\AA^{0, \le k-1}$; by lemma \ref{lem:dimsupport} it must be zero.
So $B_j \in \FF_{k-1}$, and the quotient must factor via
$E_j \onto \tau^{k-1}_\FF E_j \onto B_j$. Then $B_j$ must be the 
mdq of $\tau^{k-1}_\FF E_j$, and by induction we know that the sequence
of inclusions will terminate.

This finishes the proof of Theorem \ref{mainthm}.

\section{The space of polynomial stability conditions}
\label{sect:space}

In this section, we will describe to what extent Bridgeland's deformation
result for stability conditions carries over to polynomial stability
conditions. We will first introduce a natural topology on the space
of polynomial stability conditions (with respect to which the
stability conditions of Theorem \ref{mainthm} form a ``family'').

We will also briefly discuss what assumptions are necessary to proof a
deformation result comparable to \cite[Theorem 1.2]{Bridgeland:Stab}.

We will omit most proofs; after having adjusted all necessary
definitions, they carry over almost literally from Bridgeland's proofs.

\subsection{The topology}

We continue with the following translations of definitions of
\cite{Bridgeland:Stab} to our situation:

\begin{Def}
If the triangulated category $\DD$ is linear over a field,
a polynomial stability condition $(Z, \PP)$ on $\DD$
is called
\emph{numerical} if $Z \colon K(\DD) \to \C[m]$ factors via $\NN(\DD)$,
the numerical Grothendieck group.
\end{Def}

Let $\StabP(\DD)$ be the set of stability conditions
on $\DD$, and $\StabP^\NN(\DD)$ the subset of numerical ones.

By a \emph{semi-metric} on a set $\Sigma$ we denote a function
$d \colon \Sigma \times \Sigma \to [0, \infty]$ that satisfies the triangle
inequality and $d(x, x) = 0$, but is not necessarily finite or
non-zero for two distinct elements. Similarly, we call a function
$\|\cdot\| \colon V \to [0, \infty]$ on a vector
space a \emph{semi-norm} if it satisfies subadditivity and linearity
with respect to multiplication with scalars.

Bridgeland introduced the following semi-metric on the space of
$\R$-valued slicings:

For any $X \in \DD$ and an $\R$-valued slicing, let
$\phi^-_\PP (X)$ and $\phi^+_\PP(X)$ be the smallest and highest
phase appearing in the Harder-Narasimhan filtration of $X$ according to
\ref{def:slicing}(c), respectively. 
Then $d(\PP, \QQ) \in [0, \infty]$
is defined as
\[
d(\PP, \QQ) = \sup_{0 \neq X \in \DD}
		\left\{	\left| \phi_\PP^-(X) - \phi_\QQ^-(X)\right|,
		   	\left| \phi_\PP^+(X) - \phi_\QQ^+(X)\right| \right\}.
\]
Via the projection $\pi \colon S \to \R, \phi \mapsto \phi(+\infty)$,
we can pull back $d$ to get
a semi-metric $d_S$ on the space of $S$-valued slicings.

Following \cite[section 6]{Bridgeland:Stab}, we introduce a
semi-norm on the infinite-dimensional linear space
$\Hom(\KK(\DD), \C[m])$ for all $\sigma = (Z, \PP) \in \StabP(\DD)$:
\begin{align*}
\| \cdot \|_{\sigma} & \colon \Hom(\KK(\DD), \C[m]) \to [0, \infty]
								\nonumber\\
\| U \|_\sigma & = \sup \stv{\limsup_{m \to \infty}
			     \frac{\abs{U(E)(m)}}{\abs{Z(E)(m)}}}
			    {E \ \text{semistable in $\sigma$}}
\end{align*}

The next step is to show that \cite[Lemma 6.2]{Bridgeland:Stab} carries over:
For $0 < \epsilon < \frac 14$, and $\sigma = (Z, \PP) \in \StabP(\DD)$ define
$B_\epsilon(\sigma) \subset \StabP(\DD)$ as
\[
B_\epsilon(\sigma) = \stv{\tau = (\QQ, W)}
			  {\| W - Z\|_\sigma < \sin(\pi \epsilon) \ \text{and}\ 
			   d_S(\PP, \QQ) < \epsilon}.
\]
\begin{Lem}					\label{lem:norms-equivalent}
If $\tau = (\QQ, W) \in B_\epsilon(\sigma)$, then 
the semi-norms $\| \cdot \|_{\sigma}, \| \cdot \|_{\tau}$ of $\sigma$ and
$\tau$ are equivalent, i.e. there are constants
$k_1, k_2$ such that $k_1 \|U\|_\sigma < \|U\|_\tau < k_2 \|U\|_\sigma$
for all $U \in \Hom(\KK(\DD), \C[m])$.
\end{Lem}
The proof is identical to that of \cite[Lemma 6.2]{Bridgeland:Stab}.

On $\Hom(\KK(\DD), \C[m])$ we have the natural topology of point-wise
convergence; via the forgetful map $(Z, \PP) \mapsto Z$ we can pull
this back to get a system of open sets in $\StabP(\DD)$. Now equip
$\StabP(\DD)$ with the
topology generated, in the sense of a subbasis\footnote{A topology $T$ on
a set $S$ is
generated by a subbasis $\Pi$ of subsets of $S$ if open sets in $T$
are exactly
the (infinite) unions of finite intersections of sets in $\Pi$.},
by this system of open sets
and the sets $B_\epsilon(\sigma)$ defined above.

By the definition of the topology and Lemma \ref{lem:norms-equivalent}, 
the subspace 
\[ \stv{U \in \Hom(\KK(\DD), \C[m])}{\| U \|_\sigma < \infty} \]
is locally constant in $\StabP(\DD)$ and hence constant on a
connected component $\Sigma$, denoted by $V(\Sigma)$. It is equipped
with the topology generated by the topology of point-wise convergence
and the semi-norms $\| \cdot \|_\sigma$ for $\sigma \in \Sigma$ (which
are equivalent by lemma \ref{lem:norms-equivalent}); we have obtained:

\begin{Prop}				\label{prop:continuous}
For each connected component of $\Sigma \subset \StabP(\DD)$ there is
a topological vector space $V(\Sigma)$, which is a subspace of
$\Hom(\KK(\DD), \C[m])$, such that the forgetful map
$\Sigma \to V(\Sigma)$ given by $(Z, \PP) \mapsto Z$ is continuous.
\end{Prop}

Let $E$ be stable in some polynomial stability condition
$\sigma = (Z, \PP) \in \Sigma$. Then for any
$Z' \in V(\Sigma)$, the degree of $Z'(E)$ is bounded by
the degree of $Z(E)$. In particular, if $K(\DD)$ is finite dimensional,
then $V(\Sigma)$ is finite-dimensional. Further,
Bridgeland's space $\Stab(\DD)$ is a union
of connected components of $\StabP(\DD)$.

\begin{Prop}				\label{prop:locally-injective}
Suppose that $\sigma = (Z, \PP)$ and $\tau = (Z, \QQ)$ are
polynomial stability
conditions with identical central charge $Z$ and
$d_S(\PP, \QQ) < 1$. Then they are identical.
\end{Prop}
Again, the proof of \cite[Lemma 6.4]{Bridgeland:Stab} carries over literally.

Combining the two previous propositions, we obtain a natural continuous
and locally injective map
\[ \StabP(\DD) \supset \Sigma \to V(\Sigma) \subset
\Hom(K(\DD), \C[m]).\]
The discussion applies equally to numerical polynomial stability conditions:
for every connected component
$\Sigma \subset \StabP^\NN(\DD)$ there is a subspace
$V(\Sigma) \subset \Hom(\NN(\DD), \C[m])$ with the structure of a topological
vector space, such that the forgetful map
$(Z, \PP) \mapsto Z$ induces a locally injective continuous map
\[ \Sigma \to V(\Sigma).\]

\subsection{Deformations of a polynomial stability condition}

\begin{Def}
A polynomial stability condition $(Z, \PP)$ is called \emph{locally
finite} if there exists a real number $\epsilon > 0$ such that for all
$\phi \in S$, the quasi-abelian category
$\PP((\phi - \epsilon, \phi + \epsilon))$ is of finite length.
\end{Def}

Under this strong finiteness assumption, an analogue of Bridgeland's
deformation result can be proven:
\begin{Thm}				\label{thm:local-homeom}
Let $\sigma = (Z, \PP)$ be a locally finite polynomial stability
condition. Then there is an $\epsilon > 0$ such that if a group
homomorphism $W \colon \KK(\DD) \to \C[m]$ satisfies
$\|W - Z \|_\sigma < \sin (\pi \epsilon)$,
there is a locally finite stability condition
$\tau = (W, \QQ)$ with $d_S(\PP, \QQ) < \epsilon$.
\end{Thm}
In other words, a locally finite polynomial stability
condition in the connected component $\Sigma$ can be deformed uniquely
by deforming its central charge in the subspace
$V(\Sigma) \subset \Hom(K(\DD), \C[m])$, and the space
of locally finite polynomial stability conditions is a smooth manifold.

The theorem can be shown exactly along the lines of Bridgeland's proof. Since
we are not using the result in this paper, we omit the proof.

\bibliography{../all}                      % .bib-Datei
\bibliographystyle{alphaspecial}     % .bst-Datei

\end{document}